\documentclass[10pt]{amsart}
\usepackage{latexsym, amsmath,amssymb}
\usepackage{enumerate}
\usepackage{graphicx}

\usepackage{color}
\renewcommand{\epsilon}{\varepsilon}

\newtheorem{theorem}{Theorem}

\newtheorem{lemma}[theorem]{Lemma}
\newtheorem{corr}[theorem]{Corollary}

\newtheorem{proposition}[theorem]{Proposition}
\newtheorem{deff}[theorem]{Definition}

\newcommand{\bth}{\begin{theorem}}
	\newcommand{\ble}{\begin{lemma}}
		\newcommand{\bcor}{\begin{corr}}

			\newcommand{\bdeff}{\begin{deff}}

				\newcommand{\bprop}{\begin{proposition}}
					\newcommand{\ele}{\end{lemma}}
				\newcommand{\ecor}{\end{corr}}
			\newcommand{\edeff}{\end{deff}}
		
		\newcommand{\eprop}{\end{proposition}}

	\newcommand{\la}{\lambda}

	\newcommand{\supp}{\text{supp }}
	\renewcommand{\Pi}{\varPi}

	\renewcommand{\epsilon}{\varepsilon}
	\newcommand{\sgn}{{\text {sgn}}}

	\newcommand{\ls}{\lesssim}
	\newcommand{\gs}{\gtrsim}

	\numberwithin{equation}{section}
	
\begin{document}
		\title[Local estimates for the Hermite eigenfunctions]
		{Sharp Local $L^p$ estimates for the Hermite eigenfunctions}
		
		\author{Xing Wang and Cheng Zhang}
	\address{Department of Mathematics, Hunan University, Changsha, HN 410012, China}
	\email{xingwang@hnu.edu.cn}
	\address{Mathematical Sciences Center, Tsinghua University, Beijing 100084, China}
	\email{czhang98@tsinghua.edu.cn}


		\dedicatory{}

		\begin{abstract}We investigate the concentration of eigenfunctions for the Hermite operator $H=-\Delta+|x|^2$ in $\mathbb{R}^n$ by establishing local $L^p$ bounds over the compact sets with arbitrary dilations and translations. These new results  extend the local estimates by Thangavelu \cite{thangduke} and improve those derived from Koch-Tataru \cite{kt04}, and  explain the special phenomenon that the global $L^p$  bounds decrease in $p$ when $2\le p\le \frac{2n+6}{n+1}$. The key $L^2$-estimates  show that the local probabilities decrease away from the boundary $\{|x|=\la\}$, and then they satisfy Bohr's correspondence principle in any dimension. The proof uses the Hermite spectral projection operator represented by Mehler's formula for the Hermite-Schr\"odinger propagator $e^{-it H}$, and the strategy developed by Thangavelu \cite{thangduke} and Jeong-Lee-Ryu \cite{lee20}. We also exploit an explicit version of the stationary phase lemma and H\"ormander's $L^2$ oscillatory integral theorem. Using Koch-Tataru's strategy, we construct appropriate examples  to illustrate the possible concentrations and show the optimality of our local estimates.
			
		\end{abstract}
		
		\maketitle
		
		\section{Introduction}
		In the seminal work of Sogge \cite{sogge88}, he proved the $L^p$ eigenfunction bounds for elliptic operators on compact manifolds. They are related to a variable coefficient version of Stein-Tomas restriction theorem and a number of core problems in harmonic analysis and PDEs. Sogge's $L^p$ bounds are sharp on the sphere $S^n$, because of its periodic Hamilton flow and many highly concentrated eigenfunctions, such as Gaussian beams and zonal functions.  To investigate the concentration of eigenfunctions on manifolds, the $L^p$ bounds over geodesic balls and tubes have been studied, see Bourgain \cite{bourgain}, Blair-Sogge \cite{bsapde,bscmp,bsjems, bsjdg,bsinv}, Burq-G\'erard-Tzvetkov \cite{bgt}, Han \cite{han}, Hezari-Rivi\`ere \cite{HR}, Sogge \cite{soggeadv}, Sogge-Zelditch \cite{sz} and references therein.  Specifically, Sogge \cite{soggeadv} proved the following eigenfunction estimates over the geodesic balls $B(x,r)$ with center $x$ and radius $r$ 
		\begin{equation}\label{soggeball}\sup_{x\in M}\|e_\la\|_{L^2(B(x,r))}\le C r^\frac12\|e_\la\|_{L^2(M)},\ \la^{-1}\le r\le \text{Inj}\ M,\end{equation}
		where Inj $M$ is the injectivity radius of the compact manifold $M$.
		These estimates are saturated on the standard spheres by zonal functions, and can be improved under some dynamical or geometric assumption, such as $(M, g)$ having everywhere nonpositive
		curvature. 	Sogge \cite{soggeadv} also established the connection between local and global estimates \begin{equation}\label{ballcont}\left\|e_\lambda\right\|_{L^{\frac{2(n+1)}{n-1}}(M)} \leq C \lambda^{\frac{n-1}{2(n+1)}}\Big(r^{-\frac{n+1}{4}} \sup _{x \in M}\left\|e_\lambda\right\|_{L^2(B(x,r))}\Big)^{\frac{2}{n+1}}, \quad \lambda^{-1} \leq r \leq \text{Inj}\ M.\end{equation}
		One key ingredient in the proof of \eqref{ballcont} is the finite propagation speed of wave equations.

		 The Hermite operator $H=-\Delta+|x|^2$ in $\mathbb{R}^n$ shares some similar features with the spherical Laplacian, such as periodic Hamilton flow and many highly concentrated eigenfunctions, and the problem of obtaining $L^p$ eigenfunction bounds has received considerable interest in the context of Bochner-Riesz means \cite{thang, thang87, thangduke,kara, leeadv,chentams, chenadv,chencr,chenjga}, as well as unique continuation problems \cite{esca,ev, kt09}. To understand the nodal sets of the Hermite eigenfunctions in $\mathbb{R}^n$, the sizes of nodal sets in small balls have been studied, see  B\'erard-Helffer \cite{bh, bh2014}, Hanin-Zelditch-Zhou \cite{hzz, hzz2}, Beck-Hanin \cite{bb} and Jin \cite{jin}.   In this paper, we investigate the concentration of the Hermite eigenfunctions  in $\mathbb{R}^n$ by establishing sharp $L^p$  bounds over compact sets. Similar local estimates has already been considered by Thangavelu \cite{thangduke} and Koch-Tataru \cite{kt04}, motivated by the Bochner-Riesz conjecture and its local version, and the unique continuation problems.
		 
		 The Hermite functions are eigenfunctions of the one dimensional Hermite operator
		\[-h''_k+x^2h_k=(2k+1)h_k,\ k\in \mathbb{N}.\]
		They are 
		\[h_k(x)=e^{x^2/2}(-1)^k\frac{d^k}{dx^k}e^{-x^2},\ k\in \mathbb{N}.\]
		The Hermite functions form an
		orthonormal basis in $L^2(\mathbb{R})$ after normalization.
		
In dimension $n$, the Hermite eigenfunctions $e_\la$ satisfy
		$$(-\Delta+|x|^2)e_\la(x)=\la^2 e_\la(x),\ \  x\in \mathbb{R}^n,$$ 
		and a complete set of  eigenfunctions is given by
		\begin{equation}\label{basis}h_\alpha(x)=\prod_{j=1}^nh_{\alpha_j}(x_j)\end{equation}
		where the corresponding eigenvalue is $\la^2=2|\alpha|+n$, $\alpha\in \mathbb{N}^n$, $|\alpha|=\sum \alpha_j$. The multiplicity of the eigenvalue $2k+n$ is $\frac{(k+n-1)!}{k!(n-1)!}$. 
		\subsection{Global estimates}
\begin{figure}[h]
			\centering
			\includegraphics[width=0.8\textwidth]{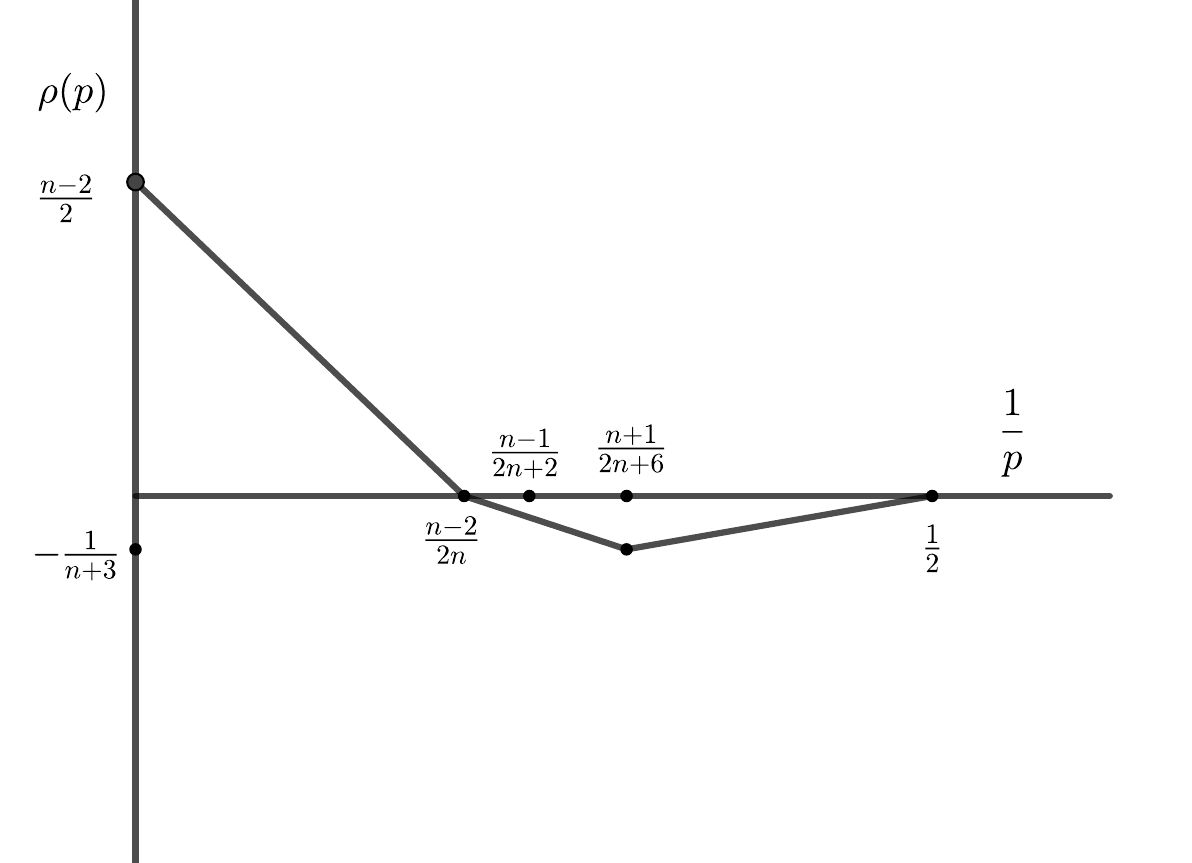}
			\caption{Global $L^p$ bounds of the Hermite eigenfunctions}
			\label{tata}
		\end{figure}

		Koch-Tataru \cite[Corollary 3.2]{kt04}  proved the following global $L^p$ eigenfunction bounds
		\begin{equation}\label{koch}	\|e_\la\|_{L^p(\mathbb{R}^n)}\ls \la^{\rho(p)}\|e_\la\|_{L^2(\mathbb{R}^n)},\end{equation}
		where for $n\ge2$,
		\begin{align}\label{rhop}
			\rho(p)=\begin{cases}
				-\frac12+\frac1p,\ \ \ \ \ 2\le p<\frac{2n+6}{n+1}  \\
				\frac{n-2}6-\frac{n}{3p},\ \ \ \frac{2n+6}{n+1}<p\le \frac{2n}{n-2}\\
				\frac{n-2}2-\frac np,\ \ \ \ \frac{2n}{n-2}< p\le \infty
			\end{cases}
		\end{align}
		and for $n=1$, 
		\begin{align}\label{rhop1}
			\rho(p)=\begin{cases}\la^{-\frac12+\frac1p},\ \ \ \ 2\le p< 4\\
				\la^{-\frac16-\frac1{3p}},\ \ \ 4< p\le \infty.\end{cases}
		\end{align}
		See Figure \ref{tata}. These results  strengthen those  of Karadzhov \cite{kara} and Thangavelu \cite{thangduke}. It is interesting to observe that the exponent $\rho(p)$ is decreasing when $p<\frac{2n+6}{n+1}$, and increasing  when $p>\frac{2n+6}{n+1}$. It is due to the special concentration features of the eigenfunctions, see the discussion after Theorem \ref{thmLp}. At the kink point $p=\frac{2n+6}{n+1}$, it is known that \cite{ktz}
		\begin{equation}\label{kink}\|e_\la\|_{L^{\frac{2n+6}{n+1}}(\mathbb{R}^n)}\ls \la^{-\frac1{n+3}}(\log\la)^{\frac{n+1}{2n+6}}\|e_\la\|_{L^2(\mathbb{R}^n)}.\end{equation}
		The log factor is necessary when $n=1$. Recently, the log loss has been removed by Jeong-Lee-Ryu \cite{lee22} for $n\ge3$. Their significant improvements are due to a new phenomenon concerning the asymmetric localization near the sphere $\la S^{n-1}$, see \cite[Theorem 1.2]{lee22}. One may also refer to \cite{lee20, leeadv, leepq} for their related works on Bochner-Riesz means and $L^p-L^q$ bounds for the Hermite spectral projection operator. 

	The Hermite eigenfunctions are essentially concentrated in the ball $\{|x|\le \la\}$ and have an exponential Airy type decay beyond this threshold. As was observed in \cite{kt04}, the behavior of eigenfunctions inside the ball $\{|x|\le \la\}$ is not very different from (a rescaling of) what happens in a bounded domain. But considerable care is required near the boundary $\{|x|=\la\}$, where the concentration scales are different. Consequently, Koch-Tataru \cite{kt04}  split the space $\mathbb{R}^n$ into overlapping dyadic parts with respect to the distance to the boundary
		\[D_j^{int}=\{\la-|x|\approx \la 2^{-2j}\},\ 1\le 2^j\le \la^\frac23\]
		\[D^{bd}=\{||x-\la|\le \la^{-\frac13}\}\]
		\[D^{ext}=\{|x|>\la+\frac12\la^{-\frac13}\}. \]
		Note that the thickness of the annulus $D_j^{int}$ is comparable to its distance to the boundary.
		
		To state the main theorem in \cite{kt04}, we define the spaces $l_\lambda^q L^p$ of functions in $\mathbb{R}^n$ with norm
		$$
		\|f\|_{l_\lambda^q L^p}^q=\|f\|_{L^p\left(D^{e x t}\right)}^q+\|f\|_{L^p\left(D^{b d}\right)}^q+\sum_{1 \leq  2^j \leq \lambda^{2 / 3}}\|f\|_{L^p\left(D_j^{i n t}\right)}^q
		$$
		with the usual modification when $q=\infty$. 
		For $x \in \mathbb{R}^n$ we let
		$$
		y=\lambda^{-\frac{2}{3}}\left(\lambda^2-|x|^2\right), \quad\langle y\rangle_{-}=1+y_{-}, \quad\langle y\rangle_{+}=1+y_{+}
		$$
	Koch-Tataru \cite{kt04} proved that	for $2\le p\le \frac{2n+2}{n-1} $, 
			\begin{equation}\label{kt1}
			\left\|\lambda^{\frac{1}{3}-\frac{n}{3}\left(\frac{1}{2}-\frac{1}{p}\right)}\langle y\rangle_{+}^{-\frac{1}{4}+\frac{n+3}{4}\left(\frac{1}{2}-\frac{1}{p}\right)}\langle y\rangle_{-}^{1-\frac{n}{2}\left(\frac{1}{2}-\frac{1}{p}\right)} e_\la\right\|_{l_\lambda^{\infty} L^p} \lesssim\|e_\la\|_{L^2},
			\end{equation}
			and for  $\frac{2n+2}{n-1} \leq p \leq \infty$,
			\begin{equation}\label{kt2}
			\left\|\lambda^{\frac{1}{3}-\frac{n}{3}\left(\frac{1}{2}-\frac{1}{p}\right)}\langle y\rangle_{+}^{\frac{1}{2}-\frac{n}{2}\left(\frac{1}{2}-\frac{1}{p}\right)}\langle y\rangle_{-}^N e_\la\right\|_{l_\lambda^{\infty} L^p} \lesssim\|e_\la\|_{L^2},\ \forall N.
			\end{equation}
		As a corollary, one can obtain  $L^p$   bounds over these dyadic annuli. For $2\le p\le \frac{2n+2}{n-1}$ 	\begin{equation}\label{shell1}\|e_\la\|_{L^p(D_j^{int})}\ls \la^{\frac1p-\frac12}2^{j(\frac{n+1}4-\frac{n+3}{2p})}\|e_\la\|_2,\ \ 1\le 2^j\le \la^\frac23,\end{equation}
		\begin{equation}\label{bdry1}\|e_\la\|_{L^p(D^{bd})}+\|e_\la\|_{L^p(D^{ext})}\ls \la^{-\frac13+\frac n3(\frac12-\frac1p)}\|e_\la\|_2\end{equation}
		and for $\frac{2n+2}{n-1}\le p\le \infty$
		\begin{equation}\label{shell2}\|e_\la\|_{L^p(D_j^{int})}\ls (\la2^{-j})^{\frac{n-2}2-\frac np}\|e_\la\|_2,\ \ 1\le 2^j\le \la^\frac23,\end{equation}
		\begin{equation}\label{bdry2}\|e_\la\|_{L^p(D^{bd})}+\|e_\la\|_{L^p(D^{ext})}\ls \la^{-\frac13+\frac n3(\frac12-\frac1p)}\|e_\la\|_2.\end{equation}
		
One may observe that the local $L^2$ estimate over the dilated ball $D_0^{int}=\{|x|<\la/4\}$ has no improvement on the trivial bound. Similarly, the local $L^{\frac{2n+6}{n+1}}$ bounds in \eqref{shell1} cannot essentially improve the global estimate \eqref{kink} over the whole space. Moreover,  the local $L^p$ bounds \eqref{shell2} over $D_0^{int}$ strengthen Thangavelu's estimates \cite[Theorem 2]{thangduke}
		\begin{equation}\label{thang}
			\|e_\la\|_{L^p(B)}\le C \la^{\frac{n-2}2-\frac np}\|e_\la\|_{L^2(\mathbb{R}^n)},\ \ \tfrac{2n+2}{n-1}\le p\le \infty.
		\end{equation}
	where $B$ is any fixed compact set in $\mathbb{R}^n$. This means that replacing $B$ by a much larger dilated ball $D_0^{int}$ does not affect the bounds. This interesting phenomenon can be explained by the existence of the  Hermite eigenfunctions with point concentration near the origin. See Section 5 and \cite[Example 5.2]{kt04}.

	\subsection{Bohr's correspondence principle}	
\begin{figure}[h]
			\centering
			\includegraphics[width=0.8\textwidth]{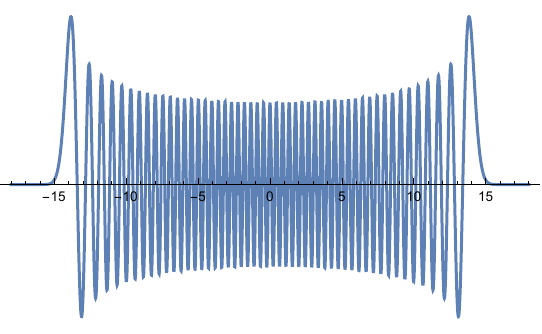}
			\caption{Hermite function $h_{100}$}
			\label{herm99}
		\end{figure}

Now let’s go back to the starting point of this study. Bohr's correspondence principle demands that classical physics and quantum physics give the same answer when the systems become large, see \cite{qc2,qc1}. For the classical harmonic oscillator, when a particle’s potential energy is equal to its energy level, it moves at its slowest speed and has the highest probability of being detected around those points. We call these points turning points, which actually correspond to the boundary $\{|x|=\la\}$. By the correspondence principle,  we would expect that for a quantum harmonic oscillator, for eigenstates with large energy, their probability density should also peak near these turning points in some sense. For the one-dimensional case, let $$s^-(x)=\int_0^x\sqrt{|t^2-\la^2|}dt,\ \ \ \  s^+(x)=\int_\la^x\sqrt{|t^2-\la^2|}dt.$$ As in Szeg\"o \cite[p. 201]{szego} and Koch-Tataru \cite[Lemma 5.1]{kt04}, the normalized eigenfunctions $\tilde h_k=h_k/\|h_k\|_2$ satisfy
		\begin{align*}
			\tilde{h}_{2 k}&=\begin{cases}
				a_{2 k}^{-}\left(\lambda^2-x^2\right)^{-\frac{1}{4}}\left(\cos s^{-}(x)+\text { error }\right) & |x|<\lambda-\lambda^{-\frac{1}{3}} \\
				O(\lambda^{-\frac{1}{6}}) & \lambda-\lambda^{-\frac{1}{3}} \leq x \leq \lambda+\lambda^{-\frac{1}{3}} \\
				a_{2 k}^{+} e^{-s^{+}(x)}\left(\lambda^2-x^2\right)^{-\frac{1}{4}}(1+\text { error }) & |x|>\lambda+\lambda^{-\frac{1}{3}}
			\end{cases} \\
			\tilde{h}_{2 k+1}&=\begin{cases}
				a_{2 k+1}^{-}\left(\lambda^2-x^2\right)^{-\frac{1}{4}}\left(\sin s^{-}(x)+\text { error }\right) & |x|<\lambda-\lambda^{-\frac{1}{3}} \\
				O(\lambda^{-\frac{1}{6}}) & \lambda-\lambda^{-\frac{1}{3}} \leq x \leq \lambda+\lambda^{-\frac{1}{3}} \\
				a_{2 k+1}^{+} e^{-s^{+}(x)}\left(\lambda^2-x^2\right)^{-\frac{1}{4}}(1+\text { error }) & x>\lambda+\lambda^{-\frac{1}{3}}
		\end{cases} \end{align*}
		where $|a_k^{ \pm}| \sim 1$, $\text { error }=O((\left|x^2-\lambda^2\right|^{-\frac{1}{2}}|| x|-\lambda|^{-1})$.
		Obviously the maximum values of Hermite eigenfunctions occur near the turning points $x=\pm\la$, see Figure \ref{herm99} by Wolfram Mathematica.
		
		 However, the eigenfunctions become more complicated in higher dimensions. Koch-Tataru \cite{kt04} showed that there are eigenfunctions that attain maximal $L^\infty$ growth near the origin. This phenomenon seems to violate Bohr's correspondence principle, and suggests that it may be not suitable to measure only by the local $L^\infty$ norm. A more reliable choice should be the local $L^2$ norm, i.e. local probability, since the square of the amplitude  is interpreted as a probability density.  As we will see in Theorem \ref{thmmain}, the local probability over  a compact set, such as the unit ball with arbitrary center, decreases as it moves away from the boundary $\{|x|=\la\}$. This satisfies Bohr's correspondence principle.
		
		\subsection{Main theorems}
		Let $B$ be any fixed compact set in $\mathbb{R}^n$, and $\nu\in \mathbb{R}^n$. For $r>0$, let
		\begin{equation}\label{ball}B(\nu,r)=\{\nu+rx:x\in B\}\end{equation}
		 be the compact set with  dilation rate $r$ and translation vector $\nu$. We shall restrict $|\nu|\le \la$ and $r\le \la$, since the Hermite eigenfunctions are essentially supported in the ball $\{|x|\le \la\}$.  The $L^p$ bounds \eqref{thang} of the Hermite eigenfunctions over the fixed compact sets  $B$ for $p\ge\frac{2n+2}{n-1}$ have been established by Thangavelu \cite{thangduke}, motivated by the local version of the  Bochner-Riesz conjecture. Koch-Tataru  \cite{kt04} strengthened those  of  Thangavelu \cite{thangduke} as well as Karadzhov \cite{kara} by obtaining global bounds over $\mathbb{R}^n$ and dyadic annuli (including the dilated ball $D_0^{int}$), and they also observed that eigenfunctions may concentrate in some  small compact subsets, e.g. balls and tubes.
		 
\begin{figure}[h]
 	\centering
 	\includegraphics[width=0.8\textwidth]{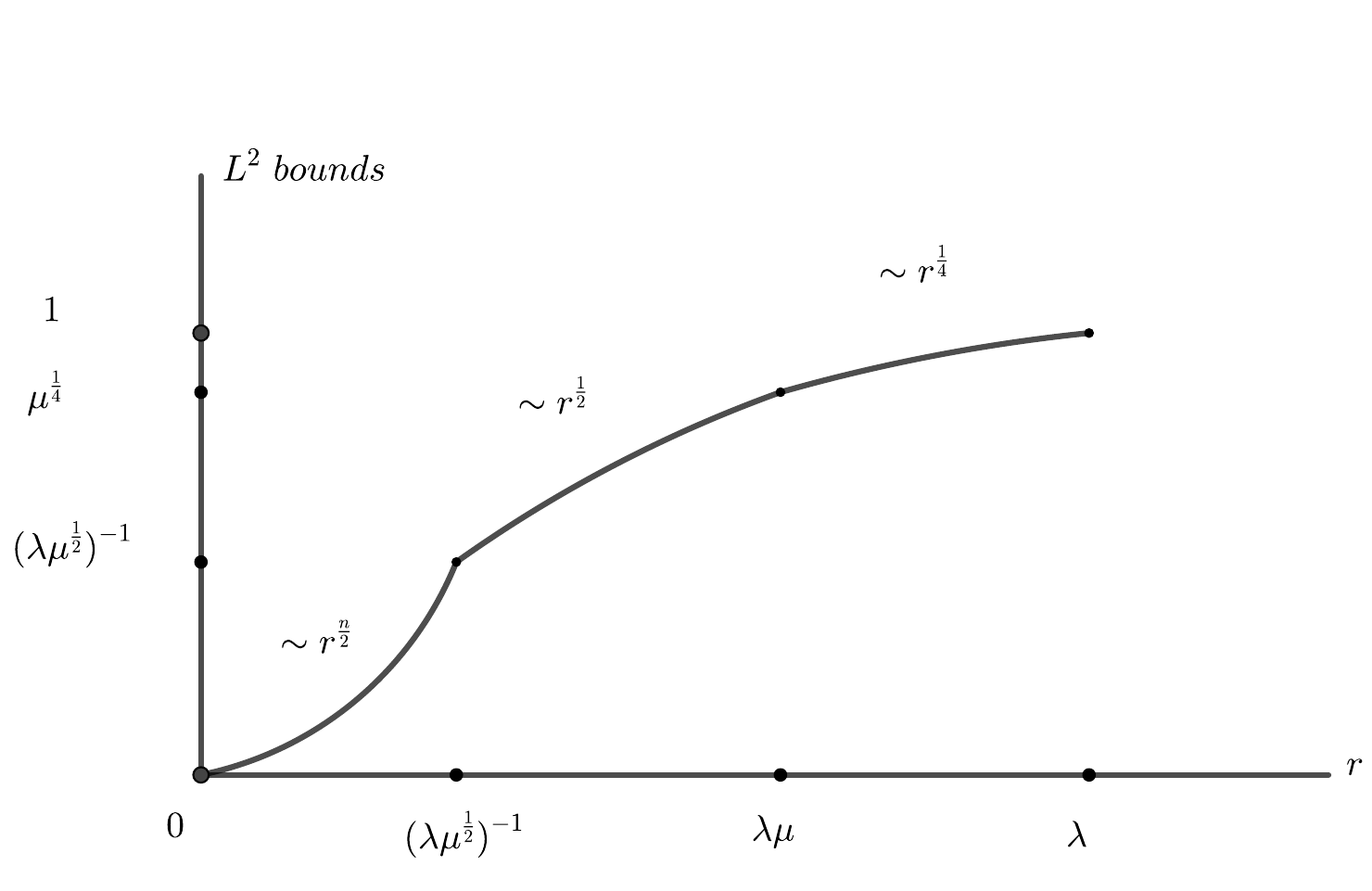}
 	\caption{Local $L^2$ bounds with respect to $r$}
 	\label{fig0}
 \end{figure}

\begin{figure}[h]
	\centering
	\includegraphics[width=0.8\textwidth]{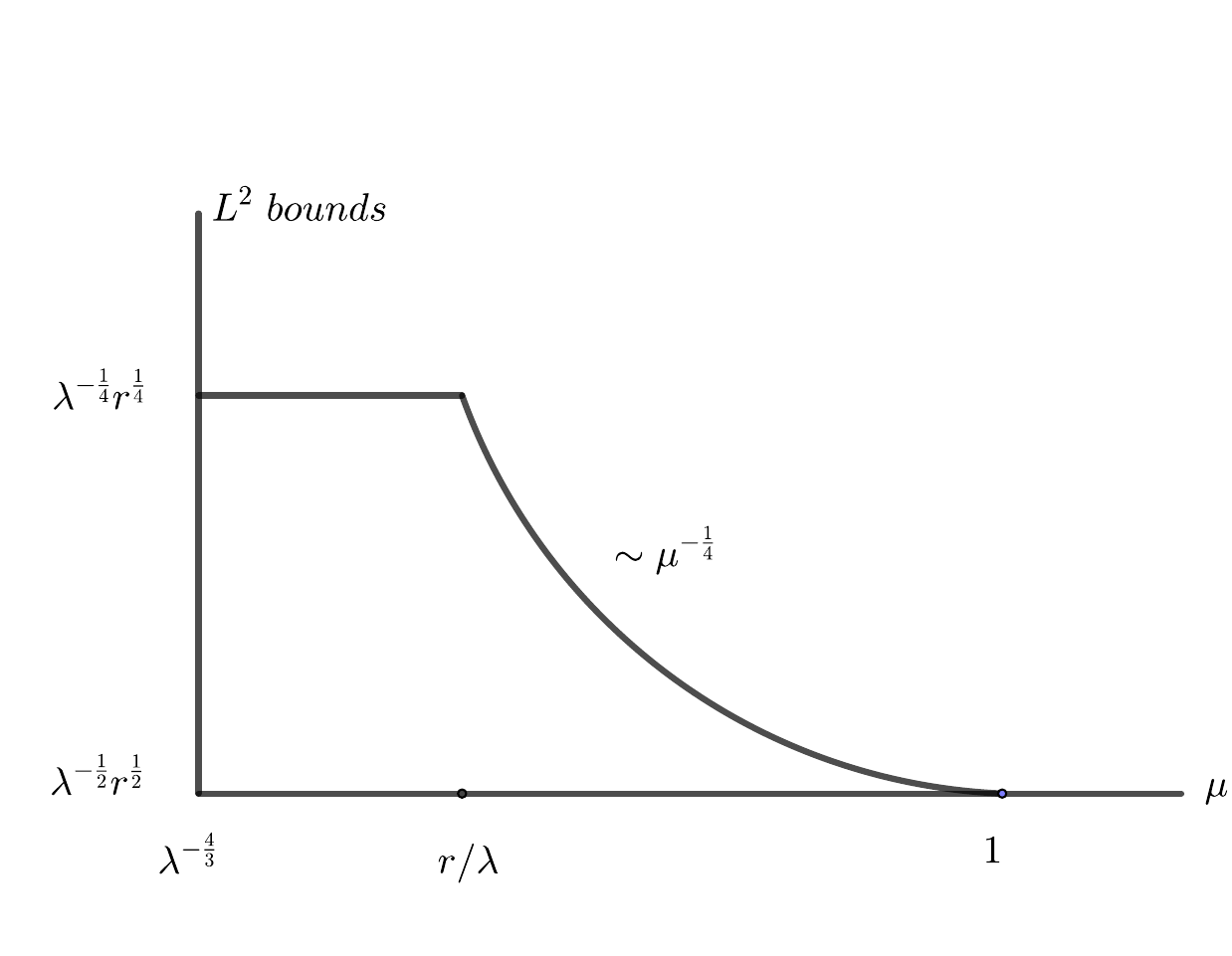}
	\caption{Local $L^2$ bounds with respect to $\mu$ when $ r\gs\la^{-\frac13}$}
	\label{fig11}
\end{figure}

\begin{figure}[h]
	\centering
	\includegraphics[width=0.8\textwidth]{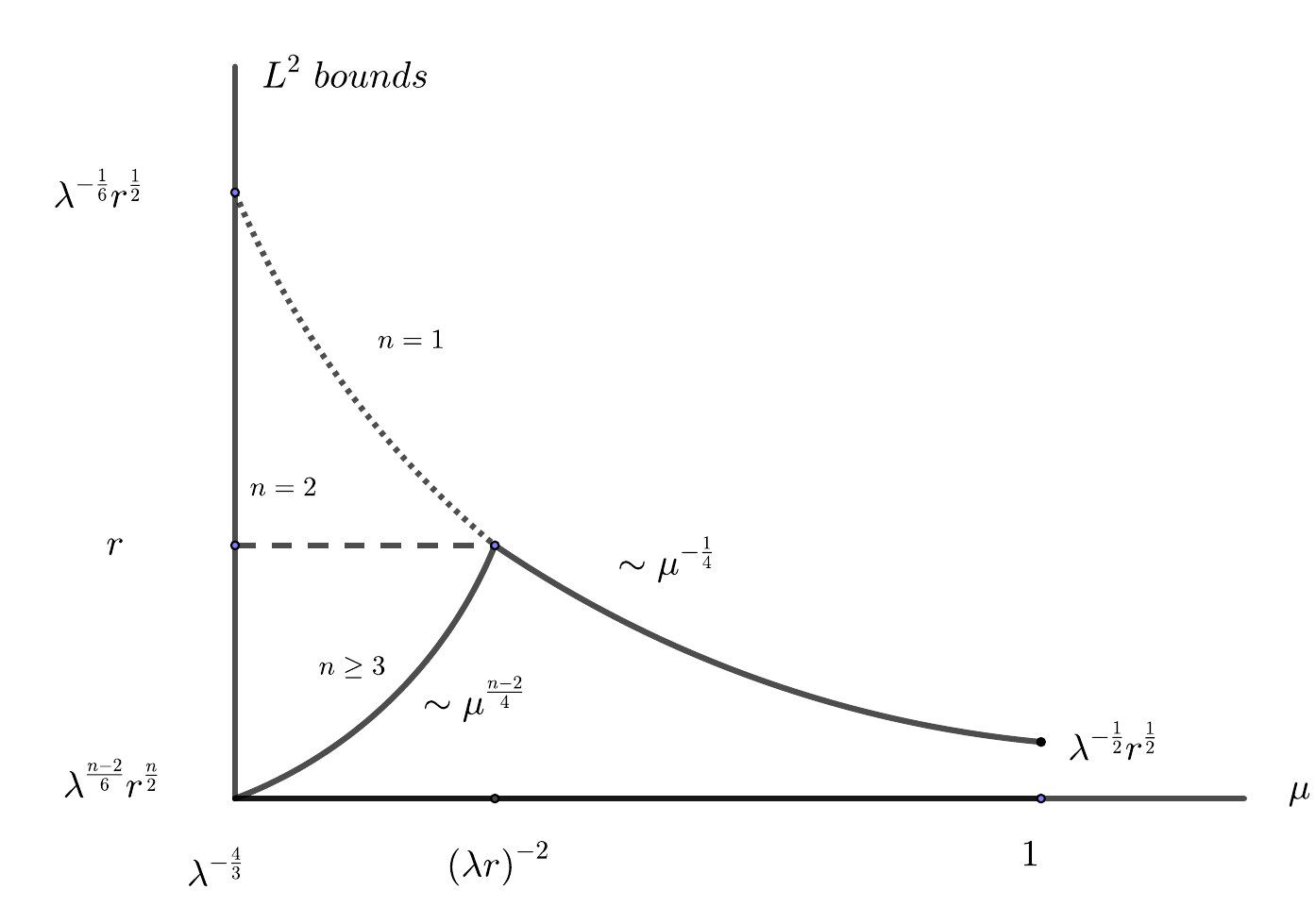}
	\caption{Local $L^2$ bounds with respect to $\mu$ when $\la^{-1}\ls r\ls\la^{-\frac13}$}
	\label{fig12}
\end{figure}

\begin{figure}[h]
	\centering
	\includegraphics[width=0.8\textwidth]{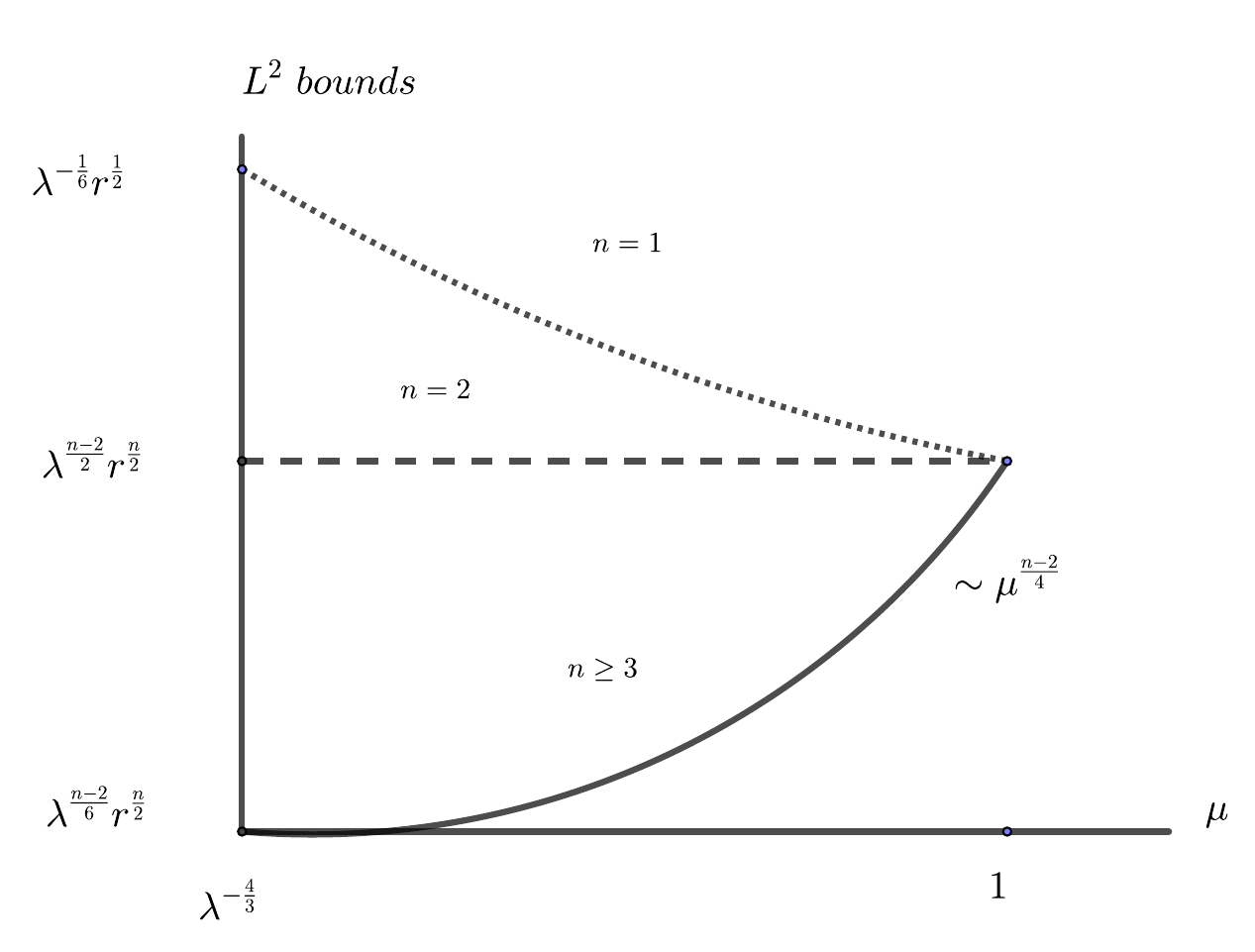}
	\caption{Local $L^2$ bounds with respect to $\mu$ when $r\ls\la^{-1}$}
	\label{fig13}
\end{figure}

		  To our best knowledge, compared to the Laplace eigenfunctions on compact manifolds, the picture of the Hermite eigenfunction estimates are still far from complete. See Section 6 for a list of related open problems. So we aim to  investigate these problems in a series of works.  In this paper, we first  establish sharp estimates over the compact sets  $B(\nu,r)$ with arbitrary dilations and translations, which extend the local estimates  of Thangavelu \cite{thangduke} and improve those derived from Koch-Tataru \cite{kt04}. Our main theorem (Theorem \ref{thmmain}) gives sharp local $L^2$ bounds (local probabilities)  over the compact sets with arbitrary  dilations and  translations. 
		 \begin{theorem}\label{thmmain}Let $B(\nu,r)$ be the compact set in \eqref{ball} and $\mu=\max\{\la^{-\frac43},1-\la^{-1}|\nu|\}$.
		 	Then for $n\ge1$ we have
		 	\begin{equation}\label{freeball}
		 		\|e_\la\|_{L^2(B(\nu,r))}\le C\Lambda(\la,r,\nu)\|e_\la\|_{L^2(\mathbb{R}^n)},
		 	\end{equation}
		 	where 
		 	\[\Lambda(\la,r,\nu)=\begin{cases}
		 		(\la\mu^\frac12)^{\frac{n-2}2}r^\frac n2,\ \ r\ls (\la\mu^\frac12)^{-1}\\
		 		(\la\mu^\frac12/r)^{-\frac12},\ \ \ (\la\mu^\frac12)^{-1}\ll r\ll \la\mu\\
		 		(\la/r)^{-\frac14},\ \ \ \ \ \ \  \  \la\mu\ls r\le \la.
		 	\end{cases}\]
		 	These bounds are sharp.\end{theorem}
 See Figures \ref{fig0}, \ref{fig11}, \ref{fig12}, \ref{fig13}.
These results  improve the bounds derived from Koch-Tataru's estimates \eqref{kt1} and \eqref{kt2} when $(\la\mu^\frac12)^{-1}\ll r\ll \la\mu$.
 The endpoints $(\la\mu^\frac12)^{-1}$ and $\la \mu$ correspond to the sizes of two different kinds of eigenfunction concentrations: point concentration and tube concentration respectively. These are similar to the eigenfunctions of the spherical Laplacian, e.g. zonal functions and Gaussian beams. Indeed, Koch-Tataru \cite{kt04} used the eigenbasis \eqref{basis} and linear combinations to construct two kinds of eigenfunctions concentrated in the dyadic annulus $D_j^{int}$. One concentrates in the ball of radius $(\la2^{-j})^{-1}$, while the other concentrates in the tube with length $\la2^{-2j}$ and radius $2^{-j/2}$. We shall exploit the strategy in \cite{kt04}  to construct new intermediate examples between these two extreme cases,  and prove the sharpness of local $L^p$ estimates, see Section 5 for the detailed construction. 	
 
 For compact sets with fixed size $r$, we shall discuss the relation between their locations and the local $L^2$ bounds. It is interesting to note that the $L^2$ bounds \eqref{freeball} with $r$ larger than the threshold $\la^{-\frac13}$ are decreasing in $\mu$ with the decay rate $\sim \mu^{-\frac14}$ when $\mu\gg r/\la$. See Figure \ref{fig11}. This suggests that the local probabilities are decreasing away from the boundary and then satisfy Bohr's correspondence principle.  When the size $r$ is smaller than the threshold $\la^{-\frac13}$, there are some significant differences between one dimension  and higher dimensions. See Figures \ref{fig12}, \ref{fig13}. It is because in the sense of $L^\infty$ norm, higher dimensional eigenfunctions can concentrate  away from the boundary, while one dimensional eigenfunctions only concentrate near the boundary.

  An important special case of Theorem \ref{thmmain} is the $\mu=r=1$ case, which extends Thangavelu's local estimates \eqref{thang} to all $p\ge2$.
 \begin{theorem}\label{coro4}Let $B$ be any fixed compact set in $\mathbb{R}^n$. Then for $n\ge1$  we have \begin{equation}\label{fixLp}\|e_\la\|_{L^p(B)}\le C\la^{\sigma(p)-\frac12}\|e_\la\|_{L^2(\mathbb{R}^n)},\end{equation}
 	where 
 	$$
 	\sigma(p)=\begin{cases}\frac{n-1}2(\frac12-\frac1p),\ \ \ 2\le p<\frac{2n+2}{n-1}\\
 		\frac{n-1}2-\frac np,\ \ \ \ \ \  \frac{2n+2}{n-1}\le p\le \infty\end{cases}
 	$$
 	is Sogge's exponent.
 	These bounds are sharp.\end{theorem}	
 Theorem \ref{coro4} is new for $2\le p<\frac{2n+2}{n-1}$. 
 It is enlightening to compare \eqref{fixLp} with Sogge’s $L^p$ estimates \cite{sogge88} on compact manifolds
 \begin{equation}\label{soggeLp}\|e_\la\|_{L^p(M)}\le C\la^{\sigma(p)}\|e_\la\|_{L^2(M)},\ 2\le p\le \infty.\end{equation}
 This phenomenon suggests that the Hermite eigenfunctions locally resemble the rescaled Laplace eigenfunctions. 
 
 Furthermore, optimal local $L^p$ estimates for all $p\ge2$ can be obtained by the interpolation between Theorem \ref{thmmain} and Koch-Tataru's $L^p$ bounds.
 \begin{theorem}\label{thmLp}
 	Let $B(\nu,r)$ be the compact set in \eqref{ball} and $\mu=\max\{\la^{-\frac43},1-\la^{-1}|\nu|\}$, and $\tilde \mu=\max\{\la^{-\frac43},\mu-\la^{-1}r\}$.
 	Then for $n\ge1$ we have
 	\begin{equation}\label{freeballLp}
 		\|e_\la\|_{L^p(B(\nu,r))}\le C\Lambda(\la,r,\nu,p)\|e_\la\|_{L^2(\mathbb{R}^n)},
 	\end{equation}
 	where for $r\ls (\la\mu^\frac12)^{-1}$
 	\[\Lambda(\la,r,\nu,p)=(\la\mu^\frac12)^{\frac{n-2}2}r^\frac np,\ \ 2\le p\le\infty,\]
 	and for $(\la\mu^\frac12)^{-1}\ll r\ll \la\mu$
 	\[\Lambda(\la,r,\nu,p)=\begin{cases}
 		(\la\mu^{\frac12}/ r)^{\frac{n-1}4-\frac{n+1}{2p}}(\la\mu^\frac12)^{\frac1p-\frac12},\ \ \ 2\le p\le \frac{2n+2}{n-1}\\
 		(\la\mu^\frac12)^{\frac{n-2}2-\frac np}\ \ \ \ \ \ \ \ \ \ \ \ \ \ \ \quad\quad \frac{2n+2}{n-1}< p\le \infty,
 	\end{cases}\]
 	and for $\la\mu\ls r\le \la$
 	\[\Lambda(\la,r,\nu,p)=\begin{cases}
 		(r/\la)^{\frac{n+3}{4p}-\frac{n+1}8}\la^{\frac1p-\frac12},\ \ 2\le p<\frac{2n+6}{n+1}\\
 		\tilde\mu^{\frac{n+3}{4p}-\frac{n+1}8}\la^{\frac1p-\frac12},\ \ \ \ \ \ \frac{2n+6}{n+1}<p\le \frac{2n+2}{n-1}\\
 		(\la\tilde \mu^\frac12)^{\frac{n-2}2-\frac np},\ \ \ \ \ \ \ \ \ \frac{2n+2}{n-1}< p\le \frac{2n}{n-2}\\
 		(\la^\frac12r^\frac12)^{\frac{n-2}2-\frac np},\ \ \ \ \ \ \ \ \frac{2n}{n-2}< p\le \infty.
 	\end{cases}\]
 	These bounds are sharp.
 \end{theorem}
 Here we use the convention that the endpoints $\frac{2n+2}{n-1}=\infty$ when $n=1$, and $\frac{2n}{n-2}=\infty$ when $n=1,2$. At the kink point $p=\frac{2n+6}{n+1}$, the local $L^p$ bounds for $\la\mu\ls r\le \la$ coincide with the global estimates \eqref{kink}. The number $\la\tilde\mu$ is essentially the distance between the compact $B(\nu,r)$ and the boundary $\{|x|=\la\}$, and obviously $\la\tilde \mu\le \la\mu$.
 
 Theorem \ref{thmLp} improves the local bounds derived from Koch-Tataru's estimates \eqref{kt1} and \eqref{kt2} when $(\la\mu^\frac12)^{-1}\ll r\ll \la\mu$ and $2\le p\le \frac{2n+2}{n-1}$. Moreover, for any fixed $r$, we can use Theorem \ref{thmLp} to determine the locations of $\nu$ (i.e. conditions on $\mu$) that  can attain the maximal local $L^p$ norms $$\sup_{\nu\in\mathbb{R}^n}\|e_\la\|_{L^p(B(\nu,r))}.$$  See the following Tables \ref{table1}, \ref{table2}, \ref{table3} with $n\ge2$. 
 
 These tables demonstrate some new concentration features of the eigenfunctions. For  $2\le p<\frac{2n+6}{n+1}$, the maximal local $L^p$ bounds significantly improve the global estimates $\la^{\frac1p-\frac12}$ whenever $r\ll \la$.  However, the maximal local $L^p$ bounds $cannot$ improve the global estimates \eqref{koch} for $p>\frac{2n+6}{n+1}$ and $r\gs \la^{-\frac13}$. This phenomenon suggests that the Hermite eigenfunctions $cannot$ highly concentrate  in any ``small'' compact set with size $r\ll\la$, in terms of the $L^p$ norm with $p$ less than $\frac{2n+6}{n+1}$. But this kind of concentration is possible for all larger $p$. These features of the eigenfunctions  explain why the exponent $\rho(p)$ in \eqref{koch} is decreasing in $p$ when $p<\frac{2n+6}{n+1}$, and increasing  when $p>\frac{2n+6}{n+1}$. See Figure \ref{tata}.
  \begin{table}[!htbp]
 	\centering
 	\caption{Maximal local estimates for $2\le p<\frac{2n+6}{n+1}$}\label{table1}	\begin{tabular}{|c|c|c|c|}
 		\hline  
 		& $\la^{-\frac13}\ls r\le \la$&$\la^{-1}\ls r\ls\la^{-\frac13}$& $r\ls \la^{-1}$\\
 		\hline  
 		max $L^p$ bounds &$(\la/r)^{\frac{n+1}8-\frac{n+3}{4p}}\la^{\frac1p-\frac12}$&$r^{\frac np-\frac{n-2}2}$&$\la^{\frac{n-2}2}r^{\frac np}$ \\
 		\hline 
 		conditions on $\mu$ & $\mu\ls r/\la$&$\mu\approx (\la r)^{-2}$&$\mu\approx 1$ \\
 		\hline
 	\end{tabular}
\end{table}
 \begin{table}[!htbp]
 	\centering
 	\caption{Maximal local estimates for $\frac{2n+6}{n+1}<p\le \frac{2n}{n-2}$}\label{table2}	\begin{tabular}{|c|c|c|c|}
 		\hline  
 		& $\la^{-\frac13}\ls r\le \la$&$\la^{-1}\ls r\ls\la^{-\frac13}$& $r\ls \la^{-1}$\\
 		\hline  
 		max $L^p$ bounds &$\la^{\frac{n-2}6-\frac n{3p}}$&$r^{\frac np-\frac{n-2}2}$&$\la^{\frac{n-2}2}r^{\frac np}$ \\
 		\hline 
 		conditions on $\mu$ & $\la\mu-r\ls \la^{-\frac13}$&$\mu\approx (\la r)^{-2}$&$\mu\approx 1$ \\
 		\hline
 	\end{tabular}
 \end{table}
 
 \begin{table}[!htbp]
 	\centering
 	\caption{Maximal local estimates for $ \frac{2n}{n-2}< p\le \infty$}\label{table3}	\begin{tabular}{|c|c|c|}
 		\hline  
 		& $\la^{-1}\ls r\le \la$& $r\ls \la^{-1}$\\
 		\hline  
 		max $L^p$ bounds &$\la^{\frac{n-2}2-\frac n{p}}$&$\la^{\frac{n-2}2}r^{\frac np}$ \\
 		\hline 
 		conditions on $\mu$ & $\mu\approx 1$&$\mu\approx 1$ \\
 		\hline
 	\end{tabular}
 \end{table}

		\subsection{Paper structure and proof sketch}
		The paper is structured as follows. In Section 2, we presents the proof of the stationary phase lemma with a precise remainder term. In Section 3, we review the representation of the kernel of the Hermite spectral projection operator. In Section 4, we prove Theorem \ref{thmmain}. In Section 5,  we show the sharpness of local $L^p$ bounds. In Section 6, we discuss some related questions for the Hermite eigenfunctions.
		
		\noindent \textbf{Proof sketch of Theorem \ref{thmmain}.}  We handle the Hermite spectral projection operator  represented by Mehler's formula for the kernel of the Hermite-Schr\"odinger propagator $e^{-itH}$. By the  strategy developed by Thangavelu \cite{thang, thang87, thangduke} and Jeong-Lee-Ryu \cite{lee22,lee20,leeadv,leepq}, we explicitly analyze  the associated oscillatory integrals by the stationary phase lemma and H\"ormander's $L^2$ oscillatory integral theorem \cite{HorL2}. The main difficulties  lie in the  discussions concerning the critical points of the phase functions in these oscillatory integrals, which require new insights. We find that in our local problem, the Hermite spectral projection operator essentially consists of two oscillatory integral operators with different phase functions, modulo some remainder terms.  One phase function behaves like the Euclidean distance function, while the other one has no obvious geometric meaning. Nevertheless, both of the phase functions satisfy the mixed Hessian condition (with rank $=n-1$) in H\"ormander's $L^2$ oscillatory integral theorem.
		
		\noindent \textbf{Proof sketch of Theorem \ref{thmLp}.} When $r\ls \la^{-1}\mu^{-\frac12}$, we may apply Koch-Tataru's $L^\infty$ estimates \eqref{shell2}, \eqref{bdry2} and then the local $L^p$ bounds follow from H\"older's inequality. When $\la\mu\ls r\le\la$, we may essentially cover the set $B(\nu,r)$ by the dyadic annuli $D_j^{int}$ with $\la\tilde\mu\ls \la2^{-2j}\ls r$ (and  the boundary  annuli, if $\tilde\mu=\la^{-\frac43}$), so the local $L^p$ bounds are implied by Koch-Tataru's $L^p$ bounds over these annuli. When $\la^{-1}\mu^{-\frac12}\ll r\ll \la\mu$, the set $B(\nu,r)$ is essentially covered by the annulus $D_j^{int}$ with $2^{-2j}\approx\mu$, so the local $L^p$ bounds follow from the interpolation between Koch-Tataru's $L^p$ bounds over $D_j^{int}$ and the local $L^2$ bounds in Theorem \ref{thmmain}.

\subsection{Notations}	Throughout this paper, $X\ls Y$ means $CX\le Y$  for some positive constants $C$ that depend only on dimension $n$ and the number of times we take derivatives and integrate by parts. In particular, if $C>1$ is a large constant and $CX\le Y$, then we denote $X\ll Y$. If $X\ls Y$ and $Y\ls X$, we denote $X\approx Y$. The notation $\|f\|_p$ means the Lebesgue norm of $f$ in $\mathbb{R}^n$. Sometimes we abbreviate the phase function $\psi(t,x,y)$ as $\psi(t)$ when $x,y$ are fixed, and  denote its partial derivative with respect to $t$ by $\psi'(t)$. 
		
		\section{Stationary phase lemmas}
		In this section, we review the one dimensional stationary phase lemma, which is important to analyze the kernel of the spectral projection operator. It is classical and there are many excellent references, e.g. H\"ormander \cite{hor}, Stein \cite{steinbook} and Sogge \cite{fio}. Using the idea in H\"ormander \cite[Theorem 7.7.5]{hor}, we prove the explicit stationary phase Lemma \ref{stphase} with a more precise remainder term than the one presented in \cite{hor}. 
				\subsection{Explicit non-stationary phase lemma}	Let  $a\in C_0^\infty(\mathbb{R})$. Let $\phi\in C^\infty(\mathbb{R})$ be real-valued. We estimate \[I_\la=\int e^{i\la\phi(t)}a(t)dt.\]
		Let $N\ge1$, $b\in C^\infty(\mathbb{R})$ and $\partial=\frac{d}{dt}$. We define
		\[(\partial\circ b)(a)=\partial(ab)=a\partial b+b\partial  a.\]
		Note that
		\[(\partial\circ b)^{N}(a)=\sum_{|\alpha|=N}c_\alpha\partial^{\alpha_0}a\partial^{\alpha_1}b\cdot\cdot\cdot\partial^{\alpha_N}b,\]
		where $\alpha=(\alpha_1,...,\alpha_N)\in \mathbb{N}^N$, and $|\alpha|=\alpha_0+...+\alpha_N$.
		For example,
		\[(\partial\circ b)^{2}(a)=(\partial\circ b)(a\partial b+b\partial a)=b^2\partial^2 a+3b\partial a\partial b+a\partial b\partial b+ab\partial^2 b.\]
		Moreover, for $j=1,...,N$,
		\[\partial^{\alpha_j} (f^{-1})=f^{-\alpha_j-1}\sum_{|\gamma^j|=\alpha_j}c_{\gamma^j} \partial^{\gamma^j_1}f\cdot\cdot\cdot\partial^{\gamma^j_{\alpha_j}}f,\]
		where $\gamma^j=(\gamma_1^j,...,\gamma^j_{\alpha_j})$.
		Here we use the convention that the sum is 1 if $\alpha_j=0$.
		
		Suppose that $\phi'\ne0$ on $\supp a$. Integrating by parts $N$ times, we have
		\begin{align*}
			I_\la&\ls \la^{-N}|\supp a|\|(\partial\circ(\phi'^{-1}))^{N}(a)\|_\infty\\
			&\ls \la^{-N}|\supp a|\sum_{|\alpha|=N}\|(\partial^{\alpha_0}a)(\phi')^{\alpha_0-2N}\prod_{j=1}^N\sum_{|\gamma^j|=\alpha_j}|\partial^{\gamma_1^j}(\phi')\cdot\cdot\cdot\partial^{\gamma_{\alpha_j}^j}(\phi')|\|_\infty\\
			&\ls \la^{-N}|\supp a|\sum_{\alpha_0=0}^N\sum_{|\beta|=N-\alpha_0}\|(\partial^{\alpha_0}a)(\phi')^{\alpha_0-2N}\partial^{\beta_1}(\phi')\cdot\cdot\cdot\partial^{\beta_{N-\alpha_0}}(\phi')\|_\infty\\
			&\ls \la^{-N}|\supp a|\sum_{\alpha_0=0}^N\sum_{\sigma}\|(\partial^{\alpha_0}a)(\phi')^{\alpha_0-2N+\sigma_0}(\phi'')^{\sigma_1}\cdot\cdot\cdot(\phi^{(N-\alpha_0+1)})^{\sigma_{N-\alpha_0}}\|_\infty
		\end{align*}
		where $\beta=(\beta_1,...,\beta_{N-\alpha_0})$, and $\sigma=(\sigma_0,...,\sigma_{N-\alpha_0})$ satisfies $\sum k\sigma_k=N-\alpha_0=\sum \sigma_k$. 
		\subsection{Explicit stationary phase lemma}
		Suppose that $\phi\in C^\infty(\mathbb{R})$ is real-valued and satisfies $\phi(0)=\phi'(0)=0,\phi''(0)\ne0$. Let \begin{equation}\label{gt}g(t)=\phi(t)-\frac12\phi''(0)t^2,\end{equation}
		\begin{equation}\label{phitheta}\phi_\theta(t)=\frac12\phi''(0)t^2+\theta g(t),\ 0\le \theta\le 1.\end{equation}
		Then $\phi_1=\phi$. Suppose that
		\[|t|\le |\phi''(0)|/\|\phi'''\|_\infty,\ \ \forall t\in \supp a.\]
		Then
		\[|g(t)/t^3|\le \frac16\|\phi'''\|_\infty,\]
		\[|g'(t)/t^2|\le \frac12\|\phi'''\|_\infty,\]
		\[|g''(t)/t|\le\|\phi'''\|_\infty, \]
		\[|g^{(k)}(t)|\le \|\phi^{(k)}\|_\infty,\ k\ge3.\]
		\[|t/\phi_\theta'(t)|=\frac1{|\phi''(0)+\theta g'(t)/t|}\le \frac1{|\phi''(0)|-\frac12|t|\|\phi'''\|_\infty}\le \frac2{|\phi''(0)|},\]
		\[|t/\phi_\theta'(t)|=\frac1{|\phi''(0)+\theta g'(t)/t|}\ge \frac1{|\phi''(0)|+\frac12|t|\|\phi'''\|_\infty}\ge \frac2{3|\phi''(0)|}.\]
		Thus
		\[|\phi_\theta'(t)|\approx |t||\phi''(0)|,\]
		\[|\phi_\theta''(t)|=|\phi''(0)+\theta g''(t)|\le2 |\phi''(0)|,\]
		\[|\phi_\theta^{(k)}(t)|\le \|\phi^{(k)}\|_\infty,\ k\ge3.\]
		Let
		\[I_\la(\theta)=\int e^{i\la\phi_\theta(t)}a(t)dt,\ \ a\in C_0^\infty(\mathbb{R}).\]
		For $m\ge1$, we need to estimate
		\[I_\la(1)=\sum_{k=0}^{2m-1}I^{(k)}_\la(0)/k!+\frac1{(2m-1)!}\int_0^1 I_\la^{(2m)}(\theta)(1-\theta)^{2m-1}d\theta.\]
		For the last term,
		\begin{align*}
			&	|I_\la^{(2m)}(\theta)|\ls \la^{2m}\Big|\int e^{i\la\phi_\theta(t)}g(t)^{2m}a(t)d t\Big|\\
			&\ls \la^{2m} \la^{-3m}|\supp a|\sum_{\alpha_0=0}^{3m}\sum_{\sigma}\|(\partial^{\alpha_0}(g^{2m}a))(\phi_\theta')^{\alpha_0-6m+\sigma_0}(\phi_\theta'')^{\sigma_1}\cdot\cdot\cdot(\phi_\theta^{(3m-\alpha_0+1)})^{\sigma_{3m-\alpha_0}}\|_\infty\\
			&\ls\la^{-m}|\supp a|\sum_{\alpha_0=0}^{3m}\sum_{\gamma,\sigma}\|(\partial^{\gamma'}a)g^{\gamma_0}(g')^{\gamma_1}\cdot\cdot\cdot(g^{(\alpha_0)})^{\gamma_{\alpha_0}}(\phi_\theta')^{\alpha_0-6m+\sigma_0}(\phi_\theta'')^{\sigma_1}\cdot\cdot\cdot(\phi_\theta^{(3m-\alpha_0+1)})^{\sigma_{3m-\alpha_0}}\|_\infty\\
			&\ls\la^{-m}|\supp a|\sum_{\alpha_0=0}^{3m}\sum_{\gamma,\sigma}\|a^{(\gamma')}\|_\infty|\phi''(0)|^{2\alpha_0-12m+3\gamma_0+2\gamma_1+\gamma_2+2\sigma_0+\sigma_1}\|\phi'''\|_\infty^{6m-\alpha_0-2\gamma_0-\gamma_1+\gamma_3-\sigma_0+\sigma_2}\\
			&\quad\quad\quad\quad\quad\quad\quad\quad\quad\quad\quad\quad\quad\quad\quad\quad\quad\quad\quad\quad\quad\quad\quad\quad\quad\quad\quad\quad\quad\quad\quad \times\prod_{k\ge4}\|\phi^{(k)}\|_\infty^{\gamma_k+\sigma_{k-1}}\\
			&\ls \la^{-m}|\supp a|(\|a\|_\infty|\phi''(0)|^{-3m}\|\phi'''\|_\infty^{2m}+......)
		\end{align*}
		where $\gamma=(\gamma_0,...,\gamma_{\alpha_0})$ satisfies $\sum\gamma_k=2m$ and $\gamma'+\sum k\gamma_k=\alpha_0$, and $\sigma=(\sigma_0,...,\sigma_{3m-\alpha_0})$ satisfies $\sum k\sigma_k=3m-\alpha_0=\sum \sigma_k$. 
		
		In particular, if $|\phi''(0)|\approx B^{-1}$, $\|\phi'''\|_\infty\approx B^{-2}$, and $\|\phi^{(k)}\|_\infty\ls B^{1-k}$ for $k\ge4$, then we have an elegant estimate
		\begin{equation}\label{last}
			|I_\la^{(2m)}(\theta)|\ls\la^{-m}|\supp a|\sum_{k=0}^{3m}\|a^{(k)}\|_\infty B^{k-m},
		\end{equation}
	by observing that
		\begin{align*}
			&-(2\alpha_0-12m+3\gamma_0+2\gamma_1+\gamma_2+2\sigma_0+\sigma_1)-2(6m-\alpha_0-2\gamma_0-\gamma_1+\gamma_3-\sigma_0+\sigma_2)\\
			&\quad\quad\quad\quad\quad\quad\quad\quad\quad\quad\quad\quad\quad\quad\quad\quad\quad\quad\quad\quad-\sum_{k\ge4}(k-1)(\gamma_k+\sigma_{k-1})\\
			&=-\sum_{k\ge0}(k-1)\gamma_k-\sum_{k\ge0}k\sigma_k\\
			&=2m-\alpha_0+\gamma'-(3m-\alpha_0)\\
			&=\gamma'-m.
		\end{align*}
		Fix $M\ge 3m-1$. For $k=0,1,...,2m-1$, let \begin{equation}\label{akt}a_k(t)=g(t)^{k}a(t),\end{equation}where $g(t)$ is defined in \eqref{gt}. Then $a_k^{(j)}(0)=0$ for $j=0,1,...,3k-1$. 
		
		Recall the Fourier transform formula for Gaussian functions (e.g. \cite[Theorem 7.6.1]{hor}) for $\omega,t\in \mathbb{R}$
		\[e^{i\frac12\omega t^2}=\frac{e^{\frac{i\pi }4\sgn \omega}}{\sqrt{2\pi|\omega|}}\int e^{-itx}e^{-i \frac{x^2}{2\omega}}dx.\]
		 Then we obtain
		\begin{align*}
			I_\la^{(k)}(0)&=(i\la)^k\int e^{i\la\frac12\phi''(0)t^2}a_k(t)dt\\
			&=(i\la)^k\frac{e^{\frac{i\pi}4\sgn \phi''(0)}}{\sqrt{2\pi|\phi''(0)|\la}}\int e^{-\frac{i}{2\phi''(0)\la}\xi^2}\hat a_k(\xi)d\xi\\
			&=(i\la)^k\frac{e^{\frac{i\pi}4\sgn \phi''(0)}}{\sqrt{2\pi|\phi''(0)|\la}}\int\Big(\sum_{j=0}^{M-1}(\frac{-i}{2\phi''(0)\la})^j\frac{\xi^{2j}}{j!}+(\frac{-i}{2\phi''(0)\la})^M\frac{\xi^{2M}}{M!}e^{i\eta}\Big)\hat a_k(\xi)d\xi,\\
			&=(i\la)^k\frac{e^{\frac{i\pi}4\sgn \phi''(0)}}{\sqrt{2\pi|\phi''(0)|\la}}\cdot (2\pi)\Big(\sum_{j=\lceil\frac{3k}2\rceil}^{M-1}(\frac{i}{2\phi''(0)\la})^j\frac{a^{(2j)}_k(0)}{j!}\Big)+ R_k
		\end{align*}
		where $\eta$ is between $0$ and $\frac{-\xi^2}{2\phi''(0)\la}$. The remainder term $R_k$ satisfies \begin{align*} 
			|R_k|&\ls\la^{k-\frac12-M} |\phi''(0)|^{-\frac12-M}\int |\xi^{2M}\hat a_k(\xi)|d\xi.
		\end{align*}
	
		In summary, we state the explicit stationary phase lemma. 
		\begin{lemma}\label{stphase}{\rm  Let  $\phi\in C^\infty(\mathbb{R})$ be real-valued and satisfy $\phi(0)=\phi'(0)=0,\ \phi''(0)\ne0$.  Suppose that $a\in C_0^\infty(\mathbb{R})$ satisfies
	$\supp a\subset [-B,B]$ for some $B>0$. If $|\phi''(0)|\approx B^{-1}$, $\|\phi'''\|_\infty\approx B^{-2}$, and $\|\phi^{(k)}\|_\infty\ls B^{1-k}$ for $k\ge4$, then for $m\ge1$ and $M\ge 3m+1$, 
		we have \begin{align*}
			&\Big|\int e^{i\la\phi(t)}a(t)dt-\sum_{k=0}^{2m-1}\sum_{j=\lceil\frac{3k}2\rceil}^{M-1}c_{kj}\la^{k-j-\frac12}|\phi''(0)|^{-j-\frac12}a_k^{(2j)}(0)\Big|\\
			&\ls \la^{-m}|\supp a|\sum_{k=0}^{3m}\|a^{(k)}\|_\infty B^{k-m}\\
			&\quad\quad\quad +(\la B^{-1})^{-\frac12-M}\sum_{k=0}^{2m-1}\la^k\int_\mathbb{R} |\xi^{2M}\hat a_k(\xi)|d\xi, \end{align*}
		where
			\[a_k(t)=(\phi(t)-\tfrac12\phi''(0)t^2)^ka(t),\] \[c_{kj}=\frac{i^{k+j}\sqrt{2\pi}}{j!k!2^j}(\sgn\phi''(0))^je^{\frac{i\pi}4\sgn \phi''(0)}.\]
			Here $L^\infty$ norms are taken over $\supp a$, and $\hat a_k(\xi)$ is the Fourier transform of $a_k(t)$ on $\mathbb{R}$. The implicit constants are independent of $\la,\phi,a, B$, and only depend on $m$ and $M$.
	}
		\end{lemma}
The remainder term in Lemma \ref{stphase} is  more precise  than the one presented in \cite[Theorem 7.7.5]{hor}. Indeed, the proof in \cite{hor} uses Sobolev inequalities to simplify the integrals involving  $\hat a_k(\xi)$ and estimates them by the sup norms of the derivatives of $a$. It does not exploit the rapid decay of the Fourier transform $\hat a_k(\xi)$, which is crucial in our applications.

		\section{The spectral projection operator}
In this section, we introduce the representation formula of the Hermite spectral projection operator. It is known that the kernel of the Hermite spectral projection operator $P_\la$ can be represented by Mehler's formula for the kernel of the Hermite-Schr\"odinger propagator $e^{-itH}$. One may refer to Jeong-Lee-Ryu \cite[Section 2.1]{lee20} for a detailed introduction. Let 
\[P_\la f=\sum_{2|\alpha|+n=\la^2}\langle f,h_\alpha \rangle  h_\alpha\]
where $\{h_\alpha\}$ is an orthonormal basis in $L^2(\mathbb{R}^n)$ and $\alpha=(\alpha_1,...,\alpha_n)\in \mathbb{N}^n$. Each $h_\alpha$ is an eigenfunction of $H$ with eigenvalue $\la^2=2|\alpha|+n$, and is a tensor product of the Hermite functions on $\mathbb{R}$.

Since $\frac1{2\pi}\int_{-\pi}^\pi e^{i\frac t2s}dt=0$ for $s\in \mathbb{R}\setminus \{0\}$, formally we have
\[P_\la f=\frac1{2\pi}\int_{-\pi}^\pi e^{i\frac t2(\la^2-H)}fdt.\]
See \cite[Section 2.1]{lee20} for a detailed proof. Recall the formula (see \cite{st}, \cite[p. 11]{thang87})
\[e^{-i\frac t2H}f=c_n(\sin t)^{-\frac n2}\int_{\mathbb{R}^n} e^{ip(t,x,y)}f(y)dy\]
where $c_n=(2\pi i)^{-\frac n2}e^{in\pi/4}$,  $p(t,x,y)=\frac{a^2\cos t-2b}{2\sin t}$ with $a^2=|x|^2+|y|^2$ and $b=x\cdot y$.
Then  the kernel of $P_\la$ can be written as
		\begin{equation}
			K(x,y)=c_n\int_{-\pi}^\pi (\sin t)^{-\frac n2}e^{i\phi_\la(t,x,y)}dt
		\end{equation}
		where 
		\[\phi_\la(t,x,y)=\frac{\la^2}2t+\frac{a^2\cos t-2b}{2\sin t}.\]
		In applications, it is more convenient to deal with the rescaled kernel with $R=\la^2$
		\begin{equation}\label{K}		K_R(x,y)=K(\la x,\la y)=2c_n\int_{-\pi/2}^{\pi/2} (\sin 2t)^{-\frac n2}e^{iR\psi(t,x,y)}dt
		\end{equation}
		where 
		\[\psi(t,x,y)=t+\frac{a^2\cos 2t-2b}{2\sin 2t}.\]
		Here the square root $z^{\frac12}=\sqrt{|z|}e^{i\frac12 \arg(z)},\ \arg(z)\in(-\pi,\pi]$.
		By inserting smooth cutoff functions, we split the integral \eqref{K} into four parts
		\begin{equation}\label{fourint}\int_{-\pi/2}^{\pi/2}=\int_0^{3\pi/8}+\int_{-3\pi/8}^0+\int_{-\pi/2}^{-\pi/8}+\int_{\pi/8}^{\pi/2}=I_0+I_0^-+I_1+I_1^-.\end{equation}
		By changing variables, we obtain
		\[I_0=\int_0^{3\pi/8}(\sin 2t)^{-\frac n2}\rho_0(t)e^{iR\psi(t,x,y)}dt,\]
		\[I_0^-=e^{-i\frac {n\pi}2}\int_0^{3\pi/8}(\sin 2t)^{-\frac n2}\rho_0(t)e^{-iR\psi(t,x,y)}dt,\]
		\[I_1=e^{-i\frac\pi2(n+R)}\int_0^{3\pi/8}(\sin 2t)^{-\frac n2}\rho_0(t)e^{iR\psi(t,x,-y)}dt,\]
		\[I_1^-=e^{i\frac {\pi}2R}\int_0^{3\pi/8}(\sin 2t)^{-\frac n2}\rho_0(t)e^{-iR\psi(t,x,-y)}dt,\]
		where  $\rho_0$ is a smooth even function supported in $[-\frac {3\pi}8,\frac {3\pi}8]$ satisfying $\rho_0\equiv1$ near 0 and $\rho_0\equiv0$ near $\pm \frac{3\pi}8$, and $\rho_0(t)+\rho_0(\frac \pi2-t)\equiv1$ for $t\in[0,\frac\pi2]$. In the following, we mainly deal with $I_0$, and the other terms can be handled similarly.

		\section{Proof of Theorem \ref{thmmain}}
		In this section, we prove Theorem \ref{thmmain}. The proof uses the strategy developed by Thangavelu \cite{thang, thang87, thangduke} and Jeong-Lee-Ryu \cite{lee22,lee20,leeadv,leepq}. They handle the oscillatory integrals by  an explicit analysis on the critical points of the phase functions. By a standard $TT^*$ argument, we only need to prove the operator bound \eqref{goal1} associated with the rescaled kernel $K_R(x,y)$. Since the kernel $K_R(x,y)$ is represented as an oscillatory integral in Section 3, we analyze the critical points of the phase function, and then split the integral properly into several parts with respect to the critical points. For the parts away from the critical points, we can handle them using integration by parts and Young's inequality. The crucial parts are  those around the critical points. We need to use the stationary phase lemma in Section 2 to calculate the kernel explicitly and then apply H\"ormander's $L^2$ oscillatory integral theorem (see e.g. \cite{HorL2}, \cite[Theorem 2.1.1]{fio}). 
		
		Let $n\ge1$, $w=\la^{-1}\nu$, $\mu=\max\{\la^{-\frac43},1-|w|\}$ and $R=\la^2$. 
		
		First, we only need to prove
		\begin{equation}\label{mid}\|e_\la\|_{L^2(B(\nu,r))}\le C \la^{-\frac12}r^\frac12\mu^{-\frac14}\|e_\la\|_2,\ \ \text{for}\ \la^{-1}\mu^{-\frac12}\ll r\ll \la\mu,
			\end{equation}
since the remaining two cases $r\ls \la^{-1}\mu^{-\frac12}$ and $\la\mu\ls r\le \la$ can be handled easily. Indeed, when $r\ls \la^{-1}\mu^{-\frac12}$, we use the $L^\infty$ bounds in \eqref{shell2} and \eqref{bdry2} to obtain the desired local $L^2$ bound $(\la\mu^\frac12)^{\frac{n-2}2}r^{\frac n2}$. When $\la\mu\ls r\le \la$, we can cover $B(\nu,r)$ by some dyadic annuli and use \eqref{shell1} and \eqref{bdry1} to get the desired bound $\la^{-\frac14}r^\frac14$.

Clearly, the condition $\la^{-1}\mu^{-\frac12}\ll r\ll \la\mu$ in \eqref{mid} implies that $\mu=1-|w|\gg \la^{-\frac43}$.

		Note that
		\[\|P_\la\|^2_{L^2(\mathbb{R}^n)\to L^{2}(B(\nu,r))}=\|P_\la^*P_\la\|_{L^{2}(B(\nu,r))\to L^2(B(\nu,r))}=\|P_\la\|_{L^{2}(B(\nu,r))\to L^2(B(\nu,r))}.\]Let $T$ be the operator associated with the rescaled kernel $K_R(x,y)$
		\begin{equation}\label{KR}		K_R(x,y)=\int_{-\pi/2}^{\pi/2} (\sin 2t)^{-\frac n2}e^{iR\psi(t,x,y)}dt
		\end{equation}Note that \[\|P_\la\|_{L^{2}(B(\nu,r))\to L^{2}(B(\nu,r))}=R^{\frac n2}\|T\|_{L^{2}(B(w,R^{-\frac12}r))\to L^2(B(w,R^{-\frac12}r))}.\]
		It suffices to prove
		\begin{equation}\label{goal1}\|T\|_{L^{2}(B(w,R^{-\frac12}r))\to L^2(B(w,R^{-\frac12}r))}\ls R^{-\frac n2}(R\mu)^{-\frac12}r
		\end{equation}
		under the assumption that $R^{-\frac12}\mu^{-\frac12}\ll r\ll R^\frac12\mu$.
		\subsection{Proof of the operator bound \eqref{goal1}}
		Let $x,y\in B(w,R^{-\frac12}r)$. Let  $D(x,y)=b^2-a^2+1$. Since $\mu=1-|w|\gg R^{-\frac12}r$, we have 
			\[1-b=2\mu-\mu^2+O(R^{-\frac12}r)\]
		\[\sqrt{D}=\sqrt{(1-b)^2-|x-y|^2}=2\mu-\mu^2+O(R^{-\frac12}r).\]
		Thus,
		$1-|x|\approx 1-|y|\approx 1-b\approx \sqrt{D}\approx \mu$.

		By \eqref{shell2} and  $\la-|\nu|=\la\mu\gg r$ we have 
		\[\|P_\la\|_{L^1(B(\nu,r))\to L^\infty(B(\nu,r))}=\|P_\la\|^2_{L^2(\mathbb{R}^n)\to L^\infty(B(\nu,r))}\ls (R\mu)^{\frac{n-2}2}.\]
		After rescaling, we still have the same $L^1-L^\infty$ bound
		\[\|T\|_{L^1(B(w,R^{-\frac12}r)\to L^\infty(B(w,R^{-\frac12}r))}\ls (R\mu)^{\frac{n-2}2}.\]
		This operator bound implies the uniform bound of the kernel for $x,y\in B(w,R^{-\frac12}r)$
		\begin{equation}\label{KB}|K_R(x,y)|\ls (R\mu)^{\frac{n-2}2}.\end{equation}
		Consider the truncated kernel
		\begin{equation}\label{Kdiag}K_R(x,y)(1-\eta(|x-y|/r_0))\end{equation}
		where $\eta\in C^\infty(\mathbb{R})$ is supported in $(1,\infty)$ and equal to 1 in $(2,\infty)$, and 
		\begin{equation}\label{r0}
			r_0=R^{-\frac12}r(R\mu r^2)^{-\frac{n-1}{2n}}.
		\end{equation}
			By Young's inequality, \eqref{KB} and $\mu\gg R^{-\frac12}r$, the operator associated with the kernel \eqref{Kdiag} has $L^{2}-L^2$ norm bounded by \begin{equation}\label{young}(R\mu)^{\frac{n-2}2}r_0^{n}= R^{-\frac n2}(R\mu)^{-\frac12}r,\end{equation}
		which is desired. So in the following we only need to handle the truncated kernel
		\[K_R(x,y)\eta(|x-y|/r_0),\]
		where $|x-y|\ge r_0$ on the support. Recall that \eqref{fourint} gives
		$$K_R(x,y)=I_0+I_0^-+I_1+I_1^-.$$
		We only handle $I_0$, and other terms are similar.
		
		To analyze the oscillatory integral $I_0$, we calculate the derivative of the phase function
			\begin{equation}\label{psi1}
			-\psi'(t)\sin^22t=a^2-1-2b\cos 2t+\cos^22t.
		\end{equation}
	By solving $\psi'(t)=0$, when $b\ge\sqrt{D}$ we get two critical points $t_1, t_2$ satisfying
		\[\cos2t_1=b+\sqrt{D}\]
		\[\cos 2t_2=b-\sqrt{D}.\]
A simple calculation gives $t_1\approx |x-y|\mu^{-\frac12}$, and  $t_2\approx \mu^{\frac12}$, see \eqref{sin1} and \eqref{sin2t2} for details. In the special case $b< \sqrt{D}$, there is exactly one critical point  $t_1$. To see this, we first notice that $b+\sqrt{D}\le 1$, since $a^2-2b=|x-y|^2\ge0$.   If $b<0$, then $|x|,|y|\ls |x-y|\le R^{-\frac12}r\ll\mu$, so we have $a^2\ll\mu ^2$ and $|b|\ll\mu^2$, which gives $|b|\ll \sqrt{D}$. So we always have $0<b+\sqrt{D}\le1$. We may assume that $b\ge \sqrt{D}$ in the following, and the special case is relatively easy and can be handled  with obvious modifications.  

A crucial observation is that $1-\mu\approx 1$. Indeed, if $1-\mu=|w|\ll1$, then $a^2\ll 1$ and $|b|\ll1$. So $|b|\ll \sqrt{D}$, which contradicts to the assumption above. Later we will use this observation in \eqref{onemu} to handle the second critical point $t_2$.

We can write
		\begin{align}\label{psi1dao}\psi'(t)&=-\frac1{\sin^22t}(\cos 2t-\cos 2t_1)(\cos2t-\cos2t_2)\nonumber\\
			&=-\frac4{\sin^22t}\sin(t+t_1)\sin(t-t_1)\sin(t+t_2)\sin(t-t_2) 
		\end{align}
		
		Next, we split the integral $I_0$ into the following five parts.
		\begin{enumerate}
			\item $0< t\ll t_1$
			\item $t\approx t_1$
			\item $t_1\ll t\ll t_2$
			\item $t\approx t_2$
			\item $t_2\ll t$
		\end{enumerate}
		
		Let $\beta\in C_0^\infty(\frac12,2)$ be a Littlewood-Paley bump function, we will use it for smooth cutoff, and for dyadic decompositions several times later.
		
		We first handle the three parts $(1),(3),(5)$ using integration by parts and Young's inequality.	
		
		\noindent	\textbf{Part (1):} For $t\ll t_1$, by \eqref{psi1dao} we have
		\begin{equation}\label{dao1}|\psi^{(k)}(t)|\approx t^{-1-k}t_1^2t_2^2\approx t^{-1-k}|x-y|^2,\ k=1,2,....\end{equation}
		We split the integral 
		\begin{equation}\label{kernel1}\eta(|x-y|/r_0)\int \eta_0(t/t_1)\rho_0(t)(\sin 2t)^{-\frac n2}e^{iR\psi}dt=\sum_{2^{-\ell}>r_0}\sum_{2^{-j}\ll t_1}K_{j\ell}(x,y)\end{equation}
		$$K_{j\ell}(x,y)=\beta(2^\ell|x-y|)\int \beta(2^j t)\eta_0(t/t_1)\rho_0(t)(\sin 2t)^{-\frac n2}e^{iR\psi}dt,$$
		where $\eta_0\in C^\infty(\mathbb{R})$ is supported in $(-\infty,\frac34)$ and equal to 1 in $(-\infty,\frac1{2})$.
		
		By \eqref{dao1}, integration by parts gives
		\[|K_{j\ell}(x,y)|\ls 2^{\frac{n-2}2j}(R2^j2^{-2\ell})^{-N},\ \forall N.\]
		Then by Young's inequality the operator $T_{j\ell}$ associated with the kernel $K_{j\ell}$ has $L^{2}-L^2$ norm bounded by
		\[2^{\frac{n-2}2j}(R2^j2^{-2\ell})^{-N}\cdot 2^{-n\ell}. \]
		Since
		$2^{-j}\ll t_1\approx 2^{-\ell} \mu^{-\frac12}$,
		we have
		\begin{align*}
			\sum_{2^{-\ell}>r_0}\sum_{2^{-j}\ll t_1}\|T_{j\ell}\|_{2\to 2}&\ls R^{-\frac n2}(R\mu)^{\frac{n-2}2-\frac n2}\cdot(R\mu^\frac12 r_0)^{-N-\frac {n-2}2+n}\\
			&=R^{-\frac n2}(R\mu)^{-\frac12}r\cdot (R\mu r^2)^{-N_1-\frac12},
		\end{align*}
		which is better than desired if $N$ is large enough, since $R\mu r^2\gs 1$. Here $N_1=\frac1{2n}(N+\tfrac {n-2}2-n)>0.$
		
		\noindent	\textbf{Part (3):} For $t_1\ll t\ll t_2$, by \eqref{psi1dao} we get
		\begin{equation}\label{dao2}|\psi^{(k)}(t)|\approx t^{-1-k}t^2t_2^2\approx t^{1-k}\mu,\ k=1,2,....\end{equation}
		We dyadically decompose the integral as
		\begin{equation}\label{kernel2}\eta(|x-y|/r_0)\int \eta_1(t)\rho_0(t)(\sin 2t)^{-\frac n2}e^{iR\psi}dt=\sum_{t_1\ll 2^{-j}\ll t_2}K_j(x,y)\end{equation}
		$$K_j(x,y)=\eta(|x-y|/r_0)\int \beta(2^j t)\eta_1(t)\rho_0(t)(\sin 2t)^{-\frac n2}e^{iR\psi}dt,$$
		where $\eta_1\in C_0^\infty(\frac32t_1,\frac3{4}t_2)$ satisfies $\eta_1\equiv 1$ in $(2t_1,\frac1{2}t_2)$.
		
		By \eqref{dao2}, integration by parts gives
		\[|K_j(x,y)|\ls 2^{\frac{n-2}2j}(1+R\mu2^{-j})^{-N},\ \forall N.\]
		Since $t_1\approx |x-y|\mu^{-\frac12}\ge r_0\mu^{-\frac12}\gg R^{-1}\mu^{-1}$, we get
		\begin{align*}\sum_{t_1\ll 2^{-j}\ll t_2}|K_j(x,y)|&\ls \sum_{t_1\ll 2^{-j}\ll t_2}2^{\frac{n-2}2j}(R\mu2^{-j})^{-N}\\
			&\ls t_1^{-\frac{n-2}2}(R\mu t_1)^{-N}\\
			&	\ls (|x-y|\mu^{-\frac12})^{-\frac{n-2}2}(1+R\mu^\frac12|x-y|)^{-N}.\end{align*}
		By Young's inequality, the operator associated with the kernel given by \eqref{kernel2} has $L^{2}-L^2$ norm bounded by
		\begin{align*}\mu^{\frac{n-2}4}\int_{\mathbb{R}^n} |x|^{-\frac{n-2}2}(1+R\mu^\frac12|x|)^{-N}dx&\ls R^{-\frac n2}(R\mu)^{\frac{n-2}2-\frac n2}\\
			&= R^{-\frac n2}(R\mu)^{-\frac12}r\cdot (R\mu r^2)^{-\frac12},\end{align*}
		this is better than desired, since $R\mu r^2\gs1$.

		\noindent	\textbf{Part (5):} For $t\gg t_2$, by \eqref{psi1dao} we have
		\begin{equation}\label{dao3}|\psi^{(k)}(t)|\approx t^{3-k},\ k=1,2,....\end{equation}
		We need to  estimate the integral
		\begin{equation}\label{kernel3}\eta(|x-y|/r_0)\int \eta_2(t/t_2)\rho_0(t)(\sin 2t)^{-\frac n2}e^{iR\psi}dt\end{equation}
		where $\eta_2\in C^\infty(\mathbb{R})$ is supported  in $(\frac32,\infty)$  and equal to 1 in $(2,\infty)$. By \eqref{dao3} and $t_2\approx \mu^\frac12$,  integration by parts gives the bound for the integral
		\[t_2^{-\frac{n-2}2}(Rt_2^3)^{-N}\approx  \mu^{-\frac{n-2}4}(R\mu^{\frac32})^{-N}.\]
		Then by Young's inequality, the operator associated with the kernel given by \eqref{kernel3} has $L^{2}-L^2$ norm bounded by
		\begin{align}\nonumber
			\mu^{-\frac{n-2}4}(R\mu^{\frac32})^{-N} (R^{-\frac12}r)^{n}&=R^{-\frac n2} (R\mu)^{-\frac12}r\cdot (R\mu^{\frac32})^{-N+\frac12}(\mu r^{-4})^{-\frac{n-1}4}\\\label{Rr61}
			&\ls R^{-\frac n2}(R\mu)^{-\frac12}r
		\end{align}  if $N$ is large enough. The last inequality uses   $\mu\gg R^{-\frac12}r$ when $Rr^6\gs1$, and uses $\mu r^{-4}\gg1$ when $Rr^6\ls1$.
		
		Next, we handle the remaining two parts (2) and (4) by the explicit stationary phase lemma and H\"ormander's $L^2$ oscillatory integral theorem. 
		
		\noindent	\textbf{Part (2): Contribution of the first critical point.}
		
			For $t\approx t_1$, by \eqref{psi1dao} we obtain
		\begin{equation}\label{psit12}\psi''(t_1)=\frac{4\sqrt{D}}{\sin2t_1}\approx \mu t_1^{-1}\end{equation}
		\begin{equation}\label{psit13}|\psi'''(t)|\approx \mu t_1^{-2}\end{equation}
		\begin{equation}\label{psit1k}|\psi^{(k)}(t)|\ls \mu t_1^{1-k},\ k=2,3,....\end{equation}
		We need to estimate the kernel given by the oscillatory integral
		\begin{equation}\label{kert1}\eta(|x-y|/r_0)\int \beta(t/t_1)\rho_0(t)(\sin 2t)^{-\frac n2}e^{iR\psi}dt.\end{equation}
		
		Let $$a(t)=\beta(t/t_1)\rho_0(t)(\sin 2t)^{-\frac n2}$$ $$g(t)=\psi(t)-\psi(t_1)-\tfrac12\psi''(t_1)(t-t_1)^2$$ $$a_k(t)=g(t)^ka(t),\ k=0,1,2,....$$Thus, by \eqref{psit12}, \eqref{psit13}, \eqref{psit1k} we have
		\[|a_k^{(j)}(t)|\ls \mu^kt_1^{-\frac n2+k-j},\ j=0,1,2,...\]
		\[|\hat a_k(\xi)|\ls \mu^kt_1^{-\frac n2+k+1}(1+t_1|\xi|)^{-N},\ \forall N.\]
		Then by the stationary phase Lemma \ref{stphase}, we have the following expansion with a remainder term estimate
		\begin{align*}
			&\Big|\int e^{iR\psi(t)}a(t)dt-\sum_{k=0}^{2m-1}\sum_{j=\lceil\frac{3k}2\rceil}^{M-1}c_{kj}R^{k-j-\frac12}e^{iR\psi(t_1)}|\psi''(t_1)|^{-\frac12-j}a_k^{(2j)}(t_1)\Big|\\
			&\ls (R\mu  t_1)^{-m}t_1^{1-\frac n2}+\sum_{k=0}^{2m-1}(R\mu t_1)^{k-\frac12-M}t_1^{1-\frac n2}\\
			&\ls t_1^{1-\frac n2}((R\mu t_1)^{-m}+(R\mu t_1)^{2m-M-\frac32})\\
			&\ls t_1^{1-\frac n2}(R\mu t_1)^{-m}
		\end{align*}
		where $c_{kj}$ are constant coefficients.
		
		Since $t_1\approx |x-y|\mu^{-\frac12}$, the kernel corresponding to the remainder term is bounded by \begin{equation}\label{err}\eta(|x-y|/r_0)(R\mu)^{\frac{n-2}2}(R\mu^{\frac12}|x-y|)^{-m-\frac n2+1}.\end{equation}
		Then by Young's inequality, the operator associated  with \eqref{err} has $L^2-L^2$ norm bounded by
		\begin{align*}&(R\mu)^{\frac{n-2}2}\int_{|x|> r_0}(R\mu^{\frac12}|x|)^{-m-\frac n2+1}dx\\
			&\ls R^{-\frac n2}(R\mu)^{\frac{n-2}2-\frac n2}\cdot(R\mu^\frac12 r_0)^{-m-\frac {n-2}2+n}\\
			&=R^{-\frac n2}(R\mu)^{-\frac12}r\cdot (R\mu r^2)^{-N-\frac12},\end{align*}
		which is better than desired if $m$ is large enough, since $R\mu r^2\gs1$. Here $N=\frac1{2n}(m+\tfrac {n-2}2-n)>0.$

		Since all the terms in the expansion share the same oscillatory factor, it suffices to handle the leading term
		\begin{equation}\label{lead}\eta(|x-y|/r_0)R^{-\frac12}D^{-\frac14}(\sin2t_1)^{-\frac {n-1}2}e^{iR\psi(t_1)},\end{equation}
		where we use \eqref{psit12}.
		
		Recall that $\mu\gg R^{-\frac12}r$, $1-|x|\approx 1-|y|\approx 1-b\approx \sqrt{D}\approx \mu$, and
		\begin{equation}\label{sqrtD}|\partial_{x,y}^\alpha \sqrt{D}|\le C_\alpha \mu^{1-\alpha}.\end{equation}
		\[\]
		\noindent \textbf{Step 1 of Part (2): Analyze the first critical point.}
		
		Recall that
		\[\psi(t,x,y)=t+\frac{a^2\cos 2t-2b}{2\sin 2t},\]
		\[\cos 2t_1=b+\sqrt{D},\]
		\begin{equation}\label{psidao01}\psi_t'(t_1(x,y),x,y)\equiv0.\end{equation}
		We have
		\begin{equation}\label{sin1}\sin^22t_1=a^2-2b^2-2b\sqrt{D}=|x-y|^2\Big(1+\frac{2b}{1-b+\sqrt{D}}\Big).\end{equation} Let 
		\[\xi(x,y)=\sqrt{1+\frac{2b}{1-b+\sqrt{D}}}.\]
		Then 
		\begin{equation}\label{sin2t1}\sin 2t_1=|x-y|\xi(x,y),\end{equation}
		where $\xi\in C^\infty$ and 
		$$\xi(x,y)\approx \mu^{-\frac12}$$
		\begin{equation}\label{xidao}|\partial_{x,y}^\alpha \xi(x,y)|\ls \mu^{-\frac12-\alpha}.\end{equation}
		Thus, we obtain
		\[t_1(x,y)\approx |x-y|\mu^{-\frac12}\]
		\begin{equation}\label{t1dao}|\partial_{x,y}^\alpha t_1(x,y)|\ls |x-y|^{1-\alpha}\mu^{\frac12},\ x\ne y.\end{equation}
		Moreover, 
		\begin{equation}\label{psit1}\psi(t_1,x,y)=|x-y|\Big(\frac{(2-a^2)t_1}{2|x-y|}+\frac{a^2(t_1-\tan t_1)}{2|x-y|}+\frac{|x-y|}{2\sin 2t_1}\Big)=|x-y|\zeta(x,y),\end{equation}
		where $\zeta\in C^\infty$ and $$\zeta(x,x)=\sqrt{1-|x|^2}\approx \sqrt\mu$$  \begin{equation}\label{zetadao}|\partial^\alpha_{x,y}\zeta(x,y)|\ls \mu^{\frac12-\alpha}\end{equation}
		by \eqref{sin2t1} and \eqref{xidao}.

		Taking partial derivatives on \eqref{psidao01} we get
		\[\partial_x[\psi_t'(t_1(x,y),x,y)]=\psi_{tt}''\partial_x t_1+\psi_{tx}''\equiv0,\]
		\[\partial_x[\psi(t_1(x,y),x,y)]=\psi_t'\partial_x t_1+\psi'_x=\psi'_x.\]
		Then
		\begin{equation}\label{pxy}\partial_y\partial_x[\psi(t_1,x,y)]=\psi_{tx}''\partial_y t_1+\psi_{xy}''=-\psi_{tt}''\partial_x t_1\partial_y t_1+\psi_{xy}''.\end{equation}
		Let $I_n$ be the $n\times n$ identity matrix. Note that\[\psi_{xy}''(t_1,x,y)=-\frac1{\sin2t_1}I_n,\]
		\begin{equation}\label{psitt1}\psi_{tt}''(t_1,x,y)=\frac{4}{\sin^32t_1}(a^2\cos2t_1-b(1+\cos^22t_1))=\frac{4\sqrt{D}}{\sin 2t_1},\end{equation}
		
		\[\partial_x t_1=\frac1{2\sin 2t_1}(\frac{x-by}{\sqrt{D}}-y)=\frac{x-y\cos 2t_1}{2\sin 2t_1\sqrt{D}}\]
		\[\partial_y t_1=\frac1{2\sin 2t_1}(\frac{y-bx}{\sqrt{D}}-x)=\frac{y-x\cos 2t_1}{2\sin 2t_1\sqrt{D}}\]
		\[\partial_xt_1\partial_yt_1=\frac{(x-y\cos2t_1)(y-x\cos 2t_1)^T}{4D\sin^22t_1}\]
		Then
		\begin{equation}\label{mix}\partial_y\partial_x[\psi(t_1,x,y)]=\frac{(x-y\cos2t_1)(x\cos 2t_1-y)^T}{\sqrt{D}\sin^32t_1}-\frac1{\sin2t_1}I_n.\end{equation}
		Note that
		\[(x\cos 2t_1-y)^T(x-y\cos 2t_1)=\sqrt D\sin^2 2t_1.\] Then by \eqref{mix} the eigenvalues of the matrix  $\partial_y\partial_x[\psi(t_1,x,y)]$ are $-\frac1{\sin 2t_1}$ (multiplicity$=n-1$) and 0.
		\[\]
		
		\noindent \textbf{Step 2 of Part (2): Apply oscillatory integral theorem.}
		
		Since $x,y\in B(w,R^{-\frac12}r)$, by \eqref{lead} we need to estimate the $L^2-L^2$ norm of the operator $T$ associated with the kernel 
		\begin{align*}
			K(x,y)=R^{-\frac12}D^{-\frac14}(\sin 2t_1)^{-\frac{n-1}2}e^{iR\psi(t_1,x,y)}\chi(R^{\frac12}r^{-1}(x-w))\chi(R^{\frac12}r^{-1}(y-w))\eta(|x-y|/r_0), 
		\end{align*}where $\chi\in C_0^\infty(\mathbb{R}^n)$. Since the phase function $\psi(t_1,x,y)$ behaves like the rescaled distance function, we may use an argument similar to Burq-G\'erard-Tzvetkov \cite[Section 6]{bgt}.
		
		We split the kernel $K$ and consider the operators $T_j$ associated with the kernel
		\[K_j(x,y)=K(x,y)\beta(2^j|x-y|),\ \ r_0<2^{-j}<R^{-\frac12}r\]
		where $r_0$ is given by \eqref{r0}.
		
		Next, we introduce a partition of unity locally finite (uniformly with respect to $j$):
		\[1=\sum_{p\in \mathbb{Z}^{n}}\chi(2^j(x-w)-p),\]
		where $\chi\in C_0^\infty(\mathbb{R}^n)$. Write 
		\[K_j(x,y)=\sum_{q,\tilde q\in \mathbb{Z}^{n}}\chi(2^j(x-w)-q)K_j(x,y)\chi(2^j(y-w)-\tilde q):=\sum_{q,\tilde q\in \mathbb{Z}^{n}} K_{j,q,\tilde q}(x,y).\]
		Note that on the support of $\chi(2^j(x-w)-p)$ we have $|2^j(x-w)-p|\ls 1$. So  
		\[|q-\tilde q|\le |2^j(x-w)-q|+|2^j(x-w)-2^j(y-w)|+|2^j(y-w)-\tilde q|\ls 1\] on the support of $K_{j,q,\tilde q}(x,y)$. Let $T_{j,q,\tilde q}$ be the operator with kernel $K_{j,q,\tilde q}(x,y)$. Then we need to handle
		\[T_j=\sum_{q,\tilde q}T_{j,q,\tilde q}.\]
		Note that the kernel of $T_{j,p,\tilde p}^*T_{j,q,\tilde q}$
		\[\int \chi(2^j(x-w)-p)\overline{K_j(x,z)}\chi(2^j(z-w)-\tilde p)\chi(2^j(z-w)-q)K_j(z,y)\chi(2^j(y-w)-\tilde q) dz\]
		vanishes when $|(p,\tilde p)-(q,\tilde q)|\gs1$. Thus, Cotlar–Stein lemma (\cite[Chapter VII]{steinbook}) implies 
		\begin{equation}\label{cotlar}\|T_{j}\|_{L^2\to L^2}\le\sup_{p,\tilde p} \sum_{q,\tilde q}\sqrt{\|T_{j,p,\tilde p}^*T_{j,q,\tilde q}\|_{L^2\to L^2}}\ls \sup_{q,\tilde q}\|T_{j,q,\tilde q}\|_{L^2\to L^2}.\end{equation}
		To estimate $\|T_{j,q,\tilde q}\|_{L^2\to L^2}$, we consider the kernel of $T_{j,q,\tilde q}$ 
		\[K_{j,q,\tilde q}(x,y)=\chi(2^j(x-w)-q)K_j(x,y)\chi(2^j(y-w)-\tilde q).\]
		Let $X=2^j(x-w)$, $Y=2^j(y-w)$. Let $\tilde T_{j,q,\tilde q}$ be the operator associated with the rescaled kernel
		\begin{align*}&\tilde K_{j,q,\tilde q}(X,Y)=\chi(X-q)K_j(2^{-j}X+w,2^{-j}Y+w)\chi(Y-\tilde q)\\
			&=\chi(X-q)\chi(Y-\tilde q)\beta(|X-Y|)K(2^{-j}X+w,2^{-j}Y+w)\\
			&=\chi(X-q)\chi(Y-\tilde q)\chi(R^{\frac12}r^{-1}2^{-j}X)\chi(R^{\frac12}r^{-1}2^{-j}Y)\beta(|X-Y|)e^{iR\psi_j(X,Y)}(\sin 2\tau_{j})^{-\frac {n-1}2}R^{-\frac12}D_j^{-\frac14}\end{align*}
		where
		\[D_j(X,Y)=D(2^{-j}X+w,2^{-j}Y+w)\] $$\tau_{j}(X,Y)=t_1(2^{-j}X+w,2^{-j}Y+w)$$
		\[\psi_{j}(X,Y)=\psi(\tau_{j},2^{-j}X+w,2^{-j}Y+w).\]
		Since $|X-Y|\approx 1$ on the support of $\beta(|X-Y|)$, by \eqref{sin2t1}, \eqref{xidao}, \eqref{psit1} and \eqref{zetadao} we have 
		\begin{equation}\label{Djlim}
			\sqrt{D_j(X,Y)}\to 1-|w|^2\approx \mu
		\end{equation}
		\begin{equation}\label{t1lim}2^j\sin( 2\tau_j(X,Y))\to |X-Y|\sqrt{1+\frac{|w|^2}{2(1-|w|^2)}}\approx |X-Y|\mu^{-\frac12}\end{equation}
		\begin{equation}\label{psit1lim}2^j\psi_j(X,Y)\to |X-Y|\sqrt{1-|w|^2}\approx |X-Y|\mu^{\frac12}\end{equation}
		in $C^\infty$ topology as $j\to \infty$.
		So for large $j$ we can write the kernel as
		\begin{align*}\tilde K_{j,q,\tilde q}(X,Y)=2^{\frac{n-1}2j}\mu^{\frac{n-3}4}R^{-\frac12}&e^{iR\psi_{j}(X,Y)}A_{j}(w,X,Y)\\
		&\cdot \chi(X-q)\chi(Y-\tilde q)\chi(R^{\frac12}r^{-1}2^{-j}X)\chi(R^{\frac12}r^{-1}2^{-j}Y)\beta(|X-Y|)\end{align*}
		where $A_{j}(w,X,Y)$ is a smooth function on the support of the kernel and satisfies
		\begin{equation}\label{ampdao}|\partial^\alpha_{X,Y} A_j(w, X,Y)|\ls 1.\end{equation} The implicit constant is independent of $j,\mu,w,q,\tilde q,R,X,Y$.

		By H\"ormander's oscillatory integral theorem, we have
		\[\|\tilde T_{j,q,\tilde q}\|_{L^2\to L^2}\ls  2^{\frac{n-1}2j}R^{-\frac12}\mu^{\frac{n-3}4}(R2^{-j} \mu^\frac12)^{-\frac{n-1}2}. \]
		By rescaling,
		\[\|T_{j,q,\tilde q}\|_{L^2\to L^2}=2^{-nj}\|\tilde T_{j,q,\tilde q}\|_{L^2\to L^2}\ls R^{-\frac n2}\mu^{-\frac12}2^{-j}.\]
		By \eqref{cotlar}, we get
		\[\|T_{j}\|_{L^2\to L^2}\ls R^{-\frac n2}\mu^{-\frac12}2^{-j}.\] 
		Then
		\[\sum_{2^{-j}<R^{-1/2}r}\|T_j\|_{2\to 2}\ls R^{-\frac n2}(R\mu)^{-\frac12}r,\] which is expected.
		\[\]		
		\noindent	\textbf{Part (4): Contribution of the second critical point.}
		
			For $t\approx t_2$, by \eqref{psi1dao} we obtain
		\begin{equation}\label{psit22}|\psi''(t_2)|=\frac{4\sqrt{D}}{\sin2t_2}\approx \mu^\frac12\end{equation}
		\begin{equation}\label{psit23}|\psi'''(t)|\approx 1\end{equation}
		\begin{equation}\label{psit2k}|\psi^{(k)}(t)|\ls \mu t_2^{1-k}\ls \mu ^{\frac{3-k}2},\ k=2,3,....\end{equation}
		We need to estimate the kernel given by the oscillatory integral
		\begin{equation}\label{kert1}\eta(|x-y|/r_0)\int \beta(t/t_2)\rho_0(t)(\sin 2t)^{-\frac n2}e^{iR\psi}dt.\end{equation}
		
		Let $$a(t)=\beta(t/t_2)\rho_0(t)(\sin 2t)^{-\frac n2}$$ $$g(t)=\psi(t)-\psi(t_2)-\tfrac12\psi''(t_2)(t-t_2)^2$$ $$a_k(t)=g(t)^ka(t),\ k=0,1,2,....$$Thus, by \eqref{psit22}, \eqref{psit23}, \eqref{psit2k} we have
		\[|a_k^{(j)}(t)|\ls \mu^kt_2^{-\frac n2+k-j},\ j=0,1,2,...\]
		\[|\hat a_k(\xi)|\ls \mu^kt_2^{-\frac n2+k+1}(1+t_2|\xi|)^{-N},\ \forall N.\]
		Then by the stationary phase Lemma \ref{stphase}, we have the following expansion with a remainder term estimate
		\begin{align*}
			&\Big|\int e^{iR\psi(t)}a(t)dt-\sum_{k=0}^{2m-1}\sum_{j=\lceil\frac{3k}2\rceil}^{M-1}c_{kj}R^{k-j-\frac12}e^{iR\psi(t_2)}|\psi''(t_2)|^{-\frac12-j}a_k^{(2j)}(t_2)\Big|\\
			&\ls (R\mu  t_2)^{-m}t_2^{1-\frac n2}+\sum_{k=0}^{2m-1}(R\mu t_2)^{k-\frac12-M}t_2^{1-\frac n2}\\
			&\ls t_2^{1-\frac n2}((R\mu t_2)^{-m}+(R\mu t_2)^{2m-M-\frac32})\\
			&\ls t_2^{1-\frac n2}(R\mu t_2)^{-m}
		\end{align*}
		where $c_{kj}$ are constant coefficients.
		
		Since $t_2\approx \mu^{\frac12}$, the kernel corresponding to the remainder term is bounded by \begin{equation}\label{err2}\mu^{\frac{2-n}4}(R\mu^{\frac32})^{-m}.\end{equation}
		Then by Young's inequality, the operator associated  with \eqref{err2} has $L^2-L^2$ norm bounded by
		\begin{align}\nonumber\mu^{\frac{2-n}4}(R\mu^{\frac32})^{-m}\int_{|x|\le R^{-1/2}r} dx&\ls R^{-\frac n2}(R\mu^\frac32)^{-m}\mu^{\frac{2-n}4}r^{n}\\\nonumber
			&=R^{-\frac n2}(R\mu)^{-\frac12}r  \cdot(R\mu^\frac32)^{-m+\frac12}(\mu r^{-4})^{-\frac{n-1}4}\\\label{Rr62}
			&\ls R^{-\frac n2}(R\mu)^{-\frac12}r  
		\end{align}
		 if $m$ is large enough.  The last inequality uses   $\mu\gg R^{-\frac12}r$ when $Rr^6\gs1$, and uses $\mu r^{-4}\gg1$ when $Rr^6\ls1$. 
		
		Since all the terms in the expansion share the same oscillatory factor, it suffices to handle  the leading term
		\begin{equation}\label{lead2}\eta(|x-y|/r_0)R^{-\frac12}D^{-\frac14}(\sin2t_2)^{-\frac {n-1}2}e^{iR\psi(t_2)},\end{equation}
		where we use \eqref{psit22}. We split \eqref{lead2} into two parts
		\[K_1(x,y)=(1-\eta(|x-y|/r_0))R^{-\frac12}D^{-\frac14}(\sin2t_2)^{-\frac {n-1}2}e^{iR\psi(t_2)}\]
		\[K(x,y)=R^{-\frac12}D^{-\frac14}(\sin2t_2)^{-\frac {n-1}2}e^{iR\psi(t_2)}.\]
		We can handle $K_1$ by Young's inequality, since it is bounded and has small support. Indeed, the operator associated with $K_1$ has $L^2-L^2$ norm bounded by
		\begin{align}\nonumber R^{-\frac12}\mu^{-\frac{n+1}4}r_0^{2n/p}&= (R\mu)^{\frac{n-2}2}(R\mu^{\frac32})^{-\frac{n-1}2}r_0^{n}\\\label{Rr63}
			&\ls (R\mu)^{\frac{n-2}2}r_0^{n}\\\nonumber
			&= R^{-\frac n2}(R\mu)^{-\frac12}r. \end{align}
	Here we use \eqref{young} and $R\mu^\frac32\gs1$. So it remains to handle $K$.
		\[\]
		\noindent \textbf{Step 1 of Part (4): Analyze the second critical point.}
		
	The calculation is similar to Step 1 of Part (2), but the behavior of the phase function associated with the second critical point are very different. So we still provide all the details.	Recall that $\mu\gg R^{-\frac12}r$, $1-|x|\approx 1-|y|\approx 1-b\approx \sqrt{D}\approx \mu$ and 
		\begin{equation}\label{sqrtD2}|\partial_{x,y}^\alpha \sqrt{D}|\le C_\alpha \mu^{1-\alpha}.\end{equation}
		Moreover,
		\[\psi(t,x,y)=t+\frac{a^2\cos 2t-2b}{2\sin 2t},\]
		\[\cos 2t_2=b-\sqrt{D},\]
		\begin{equation}\label{psidao02}\psi_t'(t_2(x,y),x,y)\equiv0.\end{equation}
		Note that
		\begin{equation}\label{sin2t2}\sin^22t_2=a^2-2b^2+2b\sqrt{D}=|x-y|^2+2b(1-b+\sqrt{D}).\end{equation}
		Thus, we get
		\[t_2(x,y)\approx \mu^{\frac12}\]
		\begin{equation}\label{t2dao}|\partial_{x,y}^\alpha t_2(x,y)|\ls \mu^{\frac12-\alpha}.\end{equation}
		We may write
		\begin{equation}\label{psit2}\psi(t_2,x,y)=\frac12(2-a^2)t_2+\frac{a^2}2(t_2-\tan t_2)+\frac{|x-y|^2}{2\sin 2t_2},\end{equation}
		then we obtain
		\begin{equation}\label{psit2dao}|\partial_{x,y}^\alpha \psi(t_2(x,y),x,y)|\ls \mu^{\frac32-\alpha}.\end{equation}
		
		Taking partial derivatives on \eqref{psidao02} we get
		\[\partial_x[\psi_t'(t_2(x,y),x,y)]=\psi_{tt}''\partial_x t_2+\psi_{tx}''\equiv0,\]
		\[\partial_x[\psi(t_2(x,y),x,y)]=\psi_t'\partial_x t_2+\psi'_x=\psi'_x.\]
		Then
		\begin{equation}\label{pxy2}\partial_y\partial_x[\psi(t_2,x,y)]=\psi_{tx}''\partial_y t_2+\psi_{xy}''=-\psi_{tt}''\partial_x t_2\partial_y t_2+\psi_{xy}''.\end{equation}
		Note that\[\psi_{xy}''(t_2,x,y)=-\frac1{\sin2t_2}I_n,\]
		\begin{equation}\label{psitt2}\psi_{tt}''(t_2,x,y)=\frac{4}{\sin^32t_2}(a^2\cos2t_2-b(1+\cos^22t_2))=-\frac{4\sqrt{D}}{\sin 2t_2},\end{equation}
		
		\[\partial_x t_2=\frac1{2\sin 2t_2}(\frac{by-x}{\sqrt{D}}-y)=\frac{y\cos 2t_2-x}{2\sin 2t_2\sqrt{D}}\]
		\[\partial_y t_2=\frac1{2\sin 2t_2}(\frac{bx-y}{\sqrt{D}}-x)=\frac{x\cos 2t_2-y}{2\sin 2t_2\sqrt{D}}\]
		\[\partial_xt_2\partial_yt_2=\frac{(y\cos2t_2-x)(x\cos 2t_2-y)^T}{4D\sin^22t_2}.\]
		
		Then
		\begin{equation}\label{mix2}\partial_y\partial_x[\psi(t_2,x,y)]=\frac{(y\cos2t_2-x)(x\cos 2t_2-y)^T}{\sqrt{D}\sin^32t_2}-\frac1{\sin2t_2}I_n.\end{equation}
		Note that
		\[(x\cos 2t_2-y)^T(y\cos 2t_2-x)=\sqrt D\sin^2 2t_2.\] Then by \eqref{mix2} the eigenvalues of the matrix  $\partial_y\partial_x[\psi(t_2,x,y)]$ are $-\frac1{\sin 2t_2}$ (multiplicity$=n-1$) and 0.
		\[\]
		\noindent \textbf{Step 2 of Part (4): Apply oscillatory integral theorem.}
		
		By rotational symmetry, we may assume that  $w=(1-\mu,0,...,0)$.  Since $x,y\in B(w, R^{-\frac12}r)\subset B(w,\mu)$, we consider the operator $T$ with the kernel 
		\begin{equation}\label{ker2}
			K(x,y)=R^{-\frac12}D^{-\frac14}(\sin 2t_2)^{-\frac{n-1}2}e^{iR\psi(t_2,x,y)}\chi((x-w)/\mu)\chi((y-w)/\mu).  
		\end{equation} 
		Here $\chi\in C_0^\infty(\mathbb{R}^n)$. Unlike the first critical point, the phase function $\psi(t_2,x,y)$ does not behave like the distance function $|x-y|$, so we need to use a different argument to handle the associated oscillatory integral operator.
		
		We first simplify the matrix in \eqref{mix2}. Note that $x=w+O(R^{-\frac12}r)$, $y=w+O(R^{-\frac12}r)$. The big oh notations are in the sense of some matrix norm. We have
		\[(y\cos2t_2-x)(x\cos2t_2-y)^T=ww^T(1-\cos2t_2)^2+O(R^{-\frac12}r),\]
		where $ww^T=diag((1-\mu)^2,0,...,0)$. Moreover,
		\[1-b=2\mu-\mu^2+O(R^{-\frac12}r)\]
		\[\sqrt{D}=\sqrt{(1-b)^2-|x-y|^2}=2\mu-\mu^2+O(R^{-\frac12}r)\]
		\begin{align}\nonumber\frac{(1-\cos2t_2)^2}{\sqrt{D}\sin^22t_2}&=\frac{(1-b+\sqrt{D})^2}{\sqrt{D}(1-(b-\sqrt{D})^2)}=\frac{1-b+\sqrt{D}}{\sqrt{D}(2-(1-b)-\sqrt{D})}\\\nonumber
			&=\frac{2(2\mu-\mu^2)+O(R^{-\frac12}r)}{(2\mu-\mu^2+O(R^{-\frac12}r))(2-2(2\mu-\mu^2)+O(R^{-\frac12}r))}\\\label{onemu}
			&=(1-\mu)^{-2}+O(\mu^{-1}R^{-\frac12}r).\end{align}
		Here we use the fact that $1-\mu\approx1$. See the observation before \eqref{psi1dao}.
		
		Thus,
		\begin{align*}
			\partial_y\partial_x[\psi(t_2(x,y),x,y)]&=\frac{(y\cos2t_2-x)(x\cos 2t_2-y)^T}{\sqrt{D}\sin^32t_2}-\frac1{\sin2t_2}I_n\\
			&=\frac1{\sin2t_2}\Big(ww^T((1-\mu)^{-2}+O(\mu^{-1}R^{-\frac12}r))-I_n\Big)\\
			&=\frac1{\sin2t_2}((1-\mu)^{-2}ww^T-I_n)+O(\mu^{-\frac32}R^{-\frac12}r)\\
			&=-\frac1{\sin2t_2}diag(0,1,...,1)+O(\mu^{-\frac32}R^{-\frac12}r).
		\end{align*}
		Let $x=(x_1,x'),\ y=(y_1,y')$. Then
		\[\partial_{y'}\partial_{x'}[\psi(t_2(x,y),x,y)]=-\frac1{\sin2t_2}I_{n-1}+O(\mu^{-\frac32}R^{-\frac12}r).\]
		This implies
		\[|\det \partial_{y'}\partial_{x'}[\psi(t_2(x,y),x,y)]|\approx (\mu^{-\frac12})^{n-1},\]
		as $\mu\gg R^{-\frac12}r$.
		
		Recall the  derivatives estimates in \eqref{sqrtD}, \eqref{t2dao}, \eqref{psit2dao}
		$$|\partial^\alpha_{x,y} \sqrt{D(x,y)}|\ls \mu^{1-\alpha}$$
		$$|\partial_{x,y}^\alpha t_2(x,y)|\ls \mu^{\frac12-\alpha}$$
		$$|\partial_{x,y}^\alpha \psi(t_2(x,y),x,y)|\ls \mu^{\frac32-\alpha}.$$

		Consider the operator $T_\mu$ with the rescaled kernel
		\begin{align*}K(\mu X,\mu Y)&=R^{-\frac12}D(\mu X,\mu Y)^{-\frac14}e^{iR\psi(t_2(\mu X,\mu Y),\mu X,\mu Y)}(\sin 2t_2(\mu X,\mu Y))^{-\frac {n-1}2}\\
			&=R^{-\frac12}\mu^{-\frac{n+1}4}e^{iR\mu^{\frac32}\Phi(\mu,X,Y)}A(\mu,X,Y)\chi(X-w/\mu)\chi( Y-w/\mu)\end{align*}
		where
		\[|\det \partial_{X'}\partial_{Y'}\Phi(\mu,X,Y)|\approx 1,\]
		\[|\partial_{X,Y}^{\alpha}\Phi(\mu,X,Y)|\ls1\]
		\[|\partial_{X,Y}^{\alpha}A(\mu,X,Y)|\ls1.\]
		These bounds follow from rescaling the derivatives estimates above, and the implicit constants are independent of $\mu,w,R,X,Y$.
		
		For fixed $X_1,Y_1$, we define
		\[T_{\mu,X_1,Y_1}g(X)=\int K(\mu X,\mu Y)g(Y')dY'.\]
		Suppose that supp $f\subset B(w/\mu,R^{-\frac12}r/\mu)\subset B(w/\mu,1)$. Then by Minkovski and H\"older inequalities,
		\begin{align*}&\|T_\mu f\|_{L^2(B( w/\mu,R^{-\frac12}r/\mu)}^2=\int_{|X_1-1+\mu|<R^{-\frac12}r} \|\int_{|Y_1-1+\mu|<R^{-\frac12}r} T_{\mu,X_1,Y_1}(f(Y_1,\cdot))dY_1\|_{L^2(\mathbb{R}^{n-1})}^2dX_1\\
			&\le\int_{|X_1-1+\mu|<R^{-\frac12}r} \Big(\int_{|Y_1-1+\mu|<R^{-\frac12}r} \|T_{\mu,X_1,Y_1}(f(Y_1,\cdot))\|_{L^2(\mathbb{R}^{n-1})}dY_1\Big)^2dX_1\\
			&\ls R^{-\frac12}r\mu^{-1}\int_{|X_1-1+\mu|<R^{-\frac12}r}\int_{|Y_1-1+\mu|<R^{-\frac12}r} \|T_{\mu,X_1,Y_1}(f(Y_1,\cdot))\|_{L^2(\mathbb{R}^{n-1})}^2dY_1dX_1\\
			&\le R^{-\frac12}r\mu^{-1}\int_{|X_1-1+\mu|<R^{-\frac12}r}\int \|T_{\mu,X_1,Y_1}\|_{L^2(\mathbb{R}^{n-1})\to L^2(\mathbb{R}^{n-1})}^2 \|f(Y_1,\cdot)\|_{L^2(\mathbb{R}^{n-1})}^2dY_1dX_1\\
			&\ls R^{-1}r^2\mu^{-2} \sup_{X_1,Y_1}\|T_{\mu,X_1,Y_1}\|_{L^2(\mathbb{R}^{n-1})\to L^2(\mathbb{R}^{n-1})}^2 \|f\|_{L^2(B( w/\mu,R^{-\frac12}r/\mu))}^2.\end{align*}

		Thus, by H\"ormander's $L^2$ oscillatory integral theorem we have
		\begin{align*}
			\|T_\mu\|_{L^2(B(w/\mu,R^{-\frac12}r/\mu))\to L^2(B(w/\mu,R^{-\frac12}r/\mu))}&\ls R^{-\frac12}r\mu^{-1}\sup_{X_1,Y_1}\|T_{\mu,X_1,Y_1}\|_{L^2(\mathbb{R}^{n-1})\to L^2(\mathbb{R}^{n-1})}\\
			&\ls R^{-\frac12}r\mu^{-1}\cdot (R\mu^{\frac32})^{-\frac{n-1}2}\cdot R^{-\frac12}\mu^{-\frac{n+1}4}\\
			&=R^{-\frac {n+1}2}\mu^{-n-\frac12}r.\end{align*}
		By rescaling,
		\begin{align*}\|T\|_{L^2(B(w,R^{-\frac12}r))\to L^2(B(w,R^{-\frac12}r))}&=\mu^n\|T_\mu\|_{L^2(B(w/\mu,R^{-\frac12}/\mu))\to L^2(B(w/\mu,R^{-\frac12}/\mu))}\\
			&\ls R^{-\frac n2}(R\mu)^{-\frac12}r\end{align*}
		which is exactly the desired bound. So we complete the proof.
		\section{Proof of the sharpness}

In this section, we construct new examples to prove the sharpness  of local $L^p$ bounds. 
			The sharpness means that for each pair of $\nu,r$, there exist eigenfunctions saturating the bound \eqref{freeballLp}.  We shall use the strategy by Koch-Tataru \cite{kt04}.
\subsection{Construction of eigenfunctions in $D_j^{int}$.}
	For $1\le 2^j\le \la^\frac23$,  $\la^{-1}2^j\ll\delta\ls 2^{-j/2}$, and fixed $x_1^*>0$ with $\la-x_1^*\approx \la2^{-2j}$, let
\[\mathcal{T}_{j,\delta}=\{x=(x_1,x')\in D_j^{int}:|x_1-x_1^*|\ll \la 2^{-j}\delta^2,\ |x'|\ll \delta\}.\]
Then we can construct the normalized eigenfunctions (dependent on $\la,\ j,\ \delta$) so that for most  $x\in\mathcal{T}_{j,\delta}$ we have
\begin{equation}\label{example}\frac{e_\la(x)}{\|e_\la\|_{L^2(\mathbb{R}^n)} }\approx \la^{-\frac12}2^{j/2}\delta^{-\frac{n-1}2}.\end{equation}		
			
To see this, we shall modify the construction of Koch-Tataru's Example 5.4 \cite{kt04}. Fix a positive integer $N$ and we consider the eigenfunctions with eigenvalue $\la^2=n+2N$. We consider the set of indices
\[I=\{\alpha\in \mathbb{N}^n:|\alpha|=N,\ \alpha_k\ \text{even}\ \text{and}\ \alpha_k\approx \delta^{-2}\ \text{for}\ 2\le k\le n\}\]		
which has size $|I|\approx \delta^{-2(n-1)}$. For some subset $J$ of $I$ of comparable size, let
\[e_\la(x)=\sum_{\alpha\in J}\prod_{k=1}^{n}\tilde h_{\alpha_k}(x_k)\]
where $\tilde h_{\alpha_k}(x_k)$ are one dimensional normalized Hermite functions. By orthogonality, 
\begin{equation}\label{L2norm}\|e_\la\|_{L^2(\mathbb{R}^n)}\approx |J|^\frac12\approx \delta^{-(n-1)}.\end{equation} 
To determine the subset $J$ such that \eqref{example} holds, we need to avoid the cancellations in the summation. We shall use the asymptotic formulas in Section 1.2. Note that for any $\alpha\in I$, we have 
\begin{equation}\label{ak}\tilde h_{\alpha_k}(x_k)\approx \delta^\frac12,\ \text{for}\ \ 2\le k\le n,\end{equation} whenever $|x_k|\ll \delta$, since
\[|s^-(x_k)|=\int_0^{|x_k|}\sqrt{|t^2-(2\alpha_k+1)|}dt\ll 1.\]
Moreover, let $u=u(\alpha_1)=\sqrt{2\alpha_1+1}$ and 
\[s_u^-(x_1)=\int_0^{x_1}\sqrt{|t^2-u^2|}dt.\]
We have   for each $u$, the value of $s_u^-(x_1)$ varies in an interval of size $\approx1$ as $x_1$ varies in an interval of size $\approx \la 2^{-j}\delta^2$, since $\delta\gg \la^{-1}2^{j}$. So $|\cos(s_u^-(x_1))|\approx 1$ for most $x_1$.

Moreover, we have \begin{equation}\label{a1}|\tilde h_{\alpha_1}(x_1)|\approx \la^{-\frac12}2^{j/2},\end{equation} since
\[\la^2-x_1^2\approx \la^22^{-2j},\ \ \la^2-u^2\approx \delta^{-2}\ll\la^{2}2^{-2j}.\]
It remains to insure that the functions $\tilde h_{\alpha_1}(x_1)$ have the same sign for all $\alpha\in J$. 
Since
\[\frac{d}{du}\frac{d}{dx_1}s_u^-(x_1)=\frac{u}{\sqrt{|x_1^2-u^2|}}\approx 2^j,\]
we may integrate this with respect to $x_1$ (from $x_1^*$ to $x_1$) and then with respect to $u$ (from $u_2$ to $u_1$):
\begin{equation}\label{diff}s_{u_1}^-(x_1)-s_{u_2}^-(x_1)=s_{u_1}(x_1^*)-s_{u_2}(x_1^*)+O(2^j(u_1-u_2) (x_1-x_1^*)).\end{equation}
The remainder term is $\ll1$ since $|x_1-x_1^*|\ll \la 2^{-j}\delta^2$ and $|u_1-u_2|\ls \la^{-1}\delta^2$.

 Fix a large constant $M\gg1$. For $k=1,...,M$, let 
\[J_k=\{\alpha\in I:s_{u}^-(x_1^*)\ \text{mod}\ 2\pi\in [\tfrac {k-1}M2\pi,\tfrac{k}M2\pi]\}.\]
Since $\cup_{k=1}^MJ_k=I$, by the pigeonhole principle there is some $J_{k_0}$ with
\[|J_{k_0}|\ge \frac1M|I|.\]
So $|J_{k_0}|\approx |I|\approx \delta^{-2(n-1)}$.
Then  we have for any two indices in $J_{k_0}$, the corresponding two phases are close modulo $2\pi$ by \eqref{diff}, namely
\[\sup_{x_1:|x_1-x_1^*|\ll \la 2^{-j}\delta^2}|s_{u_1}^-(x_1)-s_{u_2}^-(x_1)|\ll 1 \mod2\pi.\]
Together with \eqref{L2norm}, \eqref{ak} and \eqref{a1}, this implies \eqref{example}.

The examples constructed above can be viewed as intermediate cases between the two kinds of extreme examples in Koch-Tataru \cite{kt04}, which correspond to $\delta\approx \la^{-1}2^j$ (point concentration) and $\delta\approx 2^{-j/2}$ (tube concentration) respectively.  Furthermore, the construction of eigenfunctions in $D^{bd}$ is essentially the same as the construction above with $2^j=\la^\frac23$.

			\subsection{Proof of the sharpness}
		Fix $\nu\in\mathbb{R}^n$ such that $\mu=\max\{\la^{-\frac43},1-\la^{-1}|\nu|\}\approx  2^{-2j}$ and $r>0$. Let $B(\nu,r)$ be the ball $\{x\in\mathbb{R}^n:|x-\nu|<r\}$.
		
	\noindent 	\textbf{Case  1:} When $r\ls \la^{-1}\mu^{-\frac12}$, we choose $\delta=\la^{-1}\mu^{-\frac12}$ and then $B(\nu,r)$ is essentially contained in the tube $\mathcal{T}_{j,\delta}$. So we get
		\begin{align*}\frac{\|e_\la\|_{L^p(B(\nu,r))}}{\|e_\la\|_{L^2(\mathbb{R}^n)}}&\gs \la^{-\frac12}2^{j/2}\delta^{-\frac{n-1}2}\cdot r^\frac np\\
			&\approx (\la\mu^\frac12)^{\frac{n-2}2}r^\frac np,\end{align*}
		which saturates \eqref{freeballLp}.
		
		\noindent 	\textbf{Case 2:}	When $\la^{-1}\mu^{-\frac12}\ll r\ll \la\mu$, we may choose $\delta=(\la\mu^\frac12/r)^{-\frac12}$. So we have  $\delta\ll r= \la\mu^\frac12\delta^2$ and then $|B(\nu,r)\cap \mathcal{T}_{j,\delta}|\approx r\delta^{n-1}$. Thus,
		\begin{align*}\frac{\|e_\la\|_{L^p(B(\nu,r))}}{\|e_\la\|_{L^2(\mathbb{R}^n)}}&\gs \la^{-\frac12}2^{j/2}\delta^{-\frac{n-1}2}\cdot|B(\nu,r)\cap \mathcal{T}_{j,\delta}|^\frac1p \\
			&\approx \la^{-\frac12}\mu^{-\frac14}r^\frac1p\delta^{(n-1)(\frac1p-\frac12)}\\
		&=(\la^{-\frac12}r^{\frac12}\mu^{-\frac14})^{\frac{n+1}p-\frac{n-1}2}(\la\mu^\frac12)^{\frac1p-\frac12}.\end{align*}
		Moreover, we can choose $\delta=\la^{-1}\mu^{-\frac12}$, and then the tube $\mathcal{T}_{j,\delta}$ is  essentially a ball covered by $B(\nu,r)$. Since $|\mathcal{T}_{j,\delta}|\approx \delta^n$,  we get
		\begin{align*}\frac{\|e_\la\|_{L^p(B(\nu,r))}}{\|e_\la\|_{L^2(\mathbb{R}^n)}}&\gs \la^{-\frac12}2^{j/2}\delta^{-\frac{n-1}2}\cdot |\mathcal{T}_{j,\delta}|^\frac1p\\
			&\approx (\la\mu^\frac12)^{\frac{n-2}2-\frac np}.\end{align*}
	These saturate \eqref{freeballLp}. 
		
			\noindent 	\textbf{Case  3:} When $\la\mu\ls r\le \la$, we fix some integer $k$ with $\la \tilde \mu\ls \la2^{-2k}\ls r$. The ball $B(\nu,r)$ may intersect the annuli $D_k^{int}$ for these $k$. If we choose $\delta=2^{-k/2}$, then the tube $\mathcal{T}_{k,\delta}$ is essentially contained in $B(\nu,r)$ and   $|\mathcal{T}_{k,\delta}|\approx \la2^{-k} \delta^{n+1}$. So we obtain
		\begin{align*}\frac{\|e_\la\|_{L^p(B(\nu,r))}}{\|e_\la\|_{L^2(\mathbb{R}^n)}}&\gs  \la^{-\frac12}2^{k/2}\delta^{-\frac{n-1}2}\cdot |\mathcal{T}_{k,\delta}|^\frac1p \\
	&\approx \la^{\frac1p-\frac12}2^{k(\frac{n+1}4-\frac{n+3}{2p})}.\end{align*}
		 Moreover, we may choose $\delta=\la^{-1}2^{k}$. Then the tube $\mathcal{T}_{k,\delta}$ is  essentially a ball covered by $B(\nu,r)$. Since $|\mathcal{T}_{j,\delta}|\approx \delta^n$,  we get
		 \begin{align*}\frac{\|e_\la\|_{L^p(B(\nu,r))}}{\|e_\la\|_{L^2(\mathbb{R}^n)}}&\gs \la^{-\frac12}2^{k/2}\delta^{-\frac{n-1}2}\cdot |\mathcal{T}_{k,\delta}|^\frac1p\\
		 	&\approx (\la2^{-k})^{\frac{n-2}2-\frac np}.\end{align*}
		 These  saturate \eqref{freeballLp} since $\tilde \mu^{\frac12}\ls 2^{-k}\ls (r/\la)^\frac12$.	We complete the proof of the sharpness.

	\section{Further discussions}

In this section, we discuss some open problems on the Hermite eigenfunction estimates.

\begin{enumerate}[(i)]

	\item \textbf{Restriction estimates on submanifolds.} Let $\Sigma$ be a totally geodesic submanifold (not necessarily compact) in $\mathbb{R}^n$, e.g. straight lines, hyperplanes passing through the origin. The question is to establish sharp $L^p$ estimates over $\Sigma$ for the Hermite eigenfunctions. See Burq-G\'erard-Tzvetkov \cite{bgt} and Hu \cite{hu} for the restriction estimates for the Laplace eigenfunctions. See also Blair \cite{blair}, Blair-Sogge \cite{bsapde,bscmp,bsjems, bsjdg,bsinv}, Bourgain \cite{bourgain}, Bourgain-Rudnick \cite{BR}, Canzani-Galkowski-Toth \cite{cgt}, Chen-Sogge \cite{chensogge}, Chen \cite{chenxuehua}, Greenleaf-Seeger \cite{gs}, Hassell-Tacy \cite{ht15}, Hezari \cite{hez},   Huang-Zhang \cite{apde}, Sogge-Xi-Zhang \cite{sxz}, Sogge-Zelditch \cite{sz}, Tataru \cite{tataru}, Xi-Zhang \cite{xzcmp}, Wang-Zhang \cite{wzadv}, Wyman \cite{wym}, Zhang \cite{zjfa}.
	
	\item \textbf{Kakeya-Nikodym type estimates.} Let $B$ be any fixed compact set in $\mathbb{R}^n$. Let $\nu\in\mathbb{R}^n$. For real numbers $r_1\ge r_2\ge ...\ge r_n>0$, let $\textbf{r}=diag(r_1,...,r_n)$. Let 
	\[B(\nu,\textbf{r})=\{\nu+\textbf{r} x:x\in B\}.\]
	The question is to establish sharp local $L^p\ (p\ge2)$  estimates over $B(\nu,\textbf{r})$ for the Hermite eigenfunctions. The model case $r_1=...=r_n=r$ has been handled in this paper. In the case $r_1\gg r_2=...=r_n$, the set $B(\nu,\textbf{r})$ is roughly a tube, and the local $L^2$ norms are similar to the Kakeya-Nikodym bounds considered in a series of works by  Blair-Sogge \cite{bsapde,bscmp,bsjems, bsjdg,bsinv}, and closely related to the restriction estimates on straight lines.
	
		\item \textbf{Bilinear and multilinear estimates.} It is interesting to investigate the sharp bilinear and multilinear  estimates for the Hermite eigenfunctions. See Burq-G\'erard-Tzvetkov \cite{bgt04, bgtbi, bgtmu} for the bilinear and multilinear $L^2$ estimates for the Laplace eigenfunctions, and their applications to nonlinear Schr\"odinger equations. See also Guo-Han-Tacy \cite{ght} for the bilinear $L^p$ estimates of quasimodes.

\end{enumerate}

\section*{Acknowledgments}
 The authors would like to thank Dr. Xiaoyan Su and Dr. Ying Wang for  helpful discussions and comments during the research. The authors would like to thank Professor Sanghyuk Lee for helpful comments on the preprint. X.W. is partially supported by a startup grant from Hunan University. C.Z. is partially supported by NSFC Grant No.1237010173 and a startup grant from Tsinghua University. 
\section*{Declarations}
\noindent \textbf{Data availability statement.} Data sharing not applicable to this article as no datasets were generated or analyzed during the current study.

\noindent \textbf{Conflict of interests.} The authors have no relevant financial or non-financial interests to disclose.

		\bibliographystyle{plain}
		
	\end{document}